\documentclass[12pt]{article}
\usepackage{amssymb}
\parskip=\medskipamount
\newtheorem{guia}{}[section]  
\newtheorem{lemma}[guia]{Lemma} 
\newtheorem{example}[guia]{Example}
\newtheorem{definition}[guia]{Definition}  
\newtheorem{proposition}[guia]{Proposition} 
\newtheorem{corollary}[guia]{Corollary}

\newtheorem{remark}[guia]{Remark}

 1   
\font\ddpp=msbm10  scaled \magstep 1  

\def\R{\hbox{\ddpp R}}               
\def\N{\hbox{\ddpp N}}    

\def\arco#1{\kern3pt \mathop{\vbox{\ialign{##\crcr\noalign{\kern1pt}
        $\braceld\leaders\vrule\hfill\leaders\vrule\hfill\bracerd$
    \crcr\noalign{\kern1pt\nointerlineskip}
        $\hfil\displaystyle{\kern-1pt#1\kern2pt}\hfil$\crcr}}}\limits}

\newenvironment{proof}{\noindent{\bf Proof~:}}{\QED\medskip}
\def\QED{\hskip0.1em\hfill\null\ \null\nobreak\hfill
\kern3pt\lower1.8pt\vbox{\hrule\hbox   {\vrule\kern1pt\vbox{\kern1.7pt
\hbox{$\scriptstyle   QED$}\kern0.2pt}\kern1pt\vrule}\hrule}}

\begin{document}

\title{Geometric description of vakonomic and nonholonomic dynamics. Comparison of solutions}
\author{
Jorge Cort\'es
\and Manuel de Le\'on 
\and David Mart{\'\i}n de Diego
\and Sonia Mart{\'\i}nez 
}
\date{\today}
\maketitle

\bigskip

\begin{abstract}
We treat the vakonomic dynamics with general constraints within a new geometric framework which will be appropriate to study optimal control problems. We compare our formulation with Vershik-Gershkovich one in the case of linear constraints. We show how nonholonomic mechanics also admits a new geometrical description wich enables us to develop an algorithm of comparison between the solutions of both dynamics. Some examples illustrating the theory are treated.
\end{abstract}

{\bf Key words:} vakonomic dynamics, nonholonomic dynamics, optimal control, symplectic geometry

{\bf AMS subject classifications:} 70F25, 49K15, 34A26

\section{Introduction}

As is well known, the application of tools from modern differential geometry in the fields of mechanics and control theory has meant a great advance in these research areas. For example, the study of the geometrical formulation of the nonholonomic equations of motion has led to a better comprehension of locomotion generation, controllability, motion planning and trajectory tracking, raising new interesting questions in these subjects (see \cite{BlCr,BlCr2,KeMu,KoMa,Le,MuSs,Os2,OsBu,OsBu2} and references therein). On the other hand, there is now a considerable amount of papers in which optimal control problems are addressed using geometric techniques.

In this context, we present a unified geometrical formulation of the dynamics of nonholonomic and vakonomic systems. Both kind of systems have the same mathematical ``ingredients": a lagrangian function and a set of nonintegrable constraints. But the way in which the equations of motion are derived differs. In the case of vakonomic systems, the dynamics is obtained through the application of constrained variational principle \cite{arnold}. In particular, an optimal control problem can be seen as a vakonomic one. The term ``vakonomic" (``variational axiomatic kind") is inherited from Kozlov \cite{Ko}, who proposed this mechanics as an alternative set of equations of motion for a physical system under the presence of nonholonomic constraints. Nonholonomic equations of motion are deduced using d'Alambert principle when the constraints are linear or affine.

The two approaches have deserved a lot of attention in recent years (see \cite{arnold,CaFa,clr,MD,LeMaMa,LeMu,Ma,VG} and references therein). Vakonomic mechanics (also called dynamical optimization subject to nonholonomic constraints) is used in mathematical economics (growth economic theory), sub-riemannian geometry, motion of microorganisms at low Reynolds number..., as we will illustrate in Subsections 3.1, 3.2 and 3.3, while, as we have mentioned above, nonholonomic mechanics has important applications to robotics, dynamics of wheeled vehicles, motion generation...

Several authors have discussed the domains of validity of each mechanics \cite{arnold,Ko,LeMu}, and in fact, the question remains not completely closed. The solutions of these dynamics do not coincide, in general, though there are examples in which nonholonomic solutions can be seen as solutions of the constrained variational problem. In recent papers \cite{Fa,LeMu} the characterization of this situation has been studied. In \cite{LeMu} Lewis and Murray introduced the example of a ball on a rotating table and shown that the subset of solutions of the nonholonomic problem is not included in the set of vakonomic ones. In \cite{Fa} Favretti obtains conditions in some particular cases for the equivalence between the two formulations.

Our project of unifying the comparison of both dynamics from a geometrical point of view has brought us to develop new geometric frameworks for vakonomic and nonholonomic mechanics, strongly inspired in the Skinner and Rusk formulation for singular lagrangians systems \cite{SR,SR2}. By means of these approaches, we are able to compare them using an algorithm which gives rise, under appropriate conditions, to a final constraint submanifold describing the nonholonomic solutions which are also vakonomic.

The paper is structured as follows. In Section 2, we obtain the equations of motion for vakonomic mechanics, assuming an admissibility condition, which permits us to present it using the restricted lagrangian to the constraint submanifold. Let us recall that from a geometrical point of view, the lagrangian is defined on the tangent bundle $TQ$ of the configuration manifold $Q$, and $M$ represents the submanifold of $TQ$ determined by the vanishing of the nonholonomic constraints. We will treat here with arbitrary submanifolds, that is, our constraints will be nonlinear in principle. It should be also remarked that we do not consider abnormal solutions.
 
Section 3 is devoted to reformulate in geometric terms vakonomic mechanics. In this section we will use as ambient space the fibred manifold $W_0=T^*Q\times_Q M$, which is in fact a subbundle of the Whitney sum $T^{*}Q \oplus TQ$ (the phase space in Skinner and Rusk approach). Since $T^*Q$ is equipped with a canonical symplectic form we can induce a presymplectic structure $\omega$ on $T^*Q\times_Q M$. Moreover, we can consider the hamiltonian function $H_{W_0}=\langle\pi_1, \pi_2 \rangle - \pi_2^* \tilde{L}$, where $\pi_1$ and $\pi_2$ are the
canonical projections, and $\tilde{L}$ is the restriction of $L$ to $M$. Then, we prove that the equations of motion of vakonomic mechanics are equivalent to solve the presymplectic hamiltonian equation
\[
i_X\omega=dH_{W_0}
\]
Since the 2-form $\omega$ is presymplectic, a constraint algorithm must be performed in order to obtain well-defined solutions of the dynamics. If the problem is consistent,
we obtain a family of explicit solutions on the final constraint submanifold. In addition, a compatibility condition is found which characterizes when the first constraint submanifold $W_1$ is symplectic, and therefore the algorithm stabilizes at the first step. Several applications to economy, locomotion and control theory, and subriemannian geometry are given in subsequent subsections.

In Section 4, we compare our approach with that of Vershik-Gershkovich \cite{VG} for vakonomic systems with linear constraints.  We prove that both are related by a convenient presymplectomorphism, so that our approach could be consider as its generalization
for the case of nonlinear constraints.

Since we want to compare vakonomic and nonholonomic dynamics, it is necessary to construct a geometrical framework for nonholonomic mechanics using a close phase space. Indeed, in Section 5 it is proved that the nonholonomic dynamics lives on a submanifold $\tilde{M}$ of $W_{0}$. In general, we have again a presymplectic system there and a constraint algorithm is needed to obtain the dynamics on the final constraint submanifold.

In Section 6, assuming that the dynamics lives on $W_1$ and $\tilde{M}$, respectively, we can compare their solutions by means of the mapping $\Upsilon: W_1 \longrightarrow \tilde{M}$, $(\alpha, v)\longmapsto(Leg_L(v), v)$. Several illustrative examples are worked
in order to show the different behaviours. It should be remarked that recent results by Favretti \cite{Fa} are reobtained with weaker hypotheses.

\section{Variational approach to constrained mechanics}

Let $Q$ be an $n$-dimensional configuration manifold and $L : TQ
\longrightarrow \R$ an autonomous lagrangian function. If $(q^{A})$, $1\leq A\leq n$, are coordinates on $Q$, we denote by $(q^{A},\dot{q}^{A})$ the natural bundle coordinates on $TQ$ in terms of which the tangent bundle projection $\tau_Q : TQ \longrightarrow Q$ reads as $\tau_Q(q^A,\dot{q}^A)=(q^A)$.

Let us suppose that the system is subject to some constraints given by a ($2n-m$)-dimensional submanifold $M$ of $TQ$, locally defined by $\Phi^\alpha = 0$, $1 \le \alpha \le m$, where $\Phi^{\alpha}:TQ \longrightarrow \R$.

Along the paper, we will assume the admissibility condition for the submanifold $M\subseteq TQ$, that is, for all $x\in M$, we have
\[
\dim T_xM^o=\dim S^* T_xM^o \, ,
\]
where $\displaystyle{S=dq^A\otimes \frac{\partial}{\partial \dot{q}^A}}$ is the canonical vertical endomorphism (see \cite{LR}). This is equivalent to say that the rank of the matrix
\[
\frac{\partial(\Phi^1,\ldots, \Phi^m)}{\partial(\dot{q}^1,\ldots, \dot{q}^n)}
\]
is $m$ for any choice of coordinates $(q^A, \dot{q}^A)$ in $TQ$.
Consequently, by the implicit function theorem, we can locally express the constraints (reordering coordinates if necessary) as
\begin{equation}\label{constraints}
\dot{q}^{\alpha}=\Psi^{\alpha}(q^A, \dot{q}^a) \, ,
\end{equation}
where $1\leq \alpha\leq m$, $m+1\leq a\leq n$ and $1\leq A\leq n$. Then, $(q^A, \dot{q}^a)$ are local coordinates for the submanifold $M$ of $TQ$.

We denote the set of twice differentiable curves connecting two given points $x, y\in Q$ as
\[
{\cal C}^{2}(x,y) = \{ c:[0,1] \longrightarrow Q \, | \; c \; \hbox{is}
\; C^{2}, \; c(0)=x \; \hbox{and} \; c(1)=y \} \, .
\]
This set is a differentiable infinite-dimensional manifold \cite{BiSnFi}. 

Let $c$ be a curve in ${\cal C}^{2}(x,y)$. A variation of $c$ is a curve $c_s$ in ${\cal C}^2(x,y)$ such that $c_0=c$. An infinitesimal variation of $c$ is the tangent vector of a variation of $c$, that is,
\[
u(t)=\frac{dc_s(t)}{ds}\Big|_{s=0} \in T_{c(t)}Q \, .
\]
The tangent space of ${\cal C}^{2}(x,y)$ at $c$ is then given by
\[
T_{c}\,{\cal C}^{2}(x,y) = \{ u : [0,1] \longrightarrow TQ \, / \; u \;
\hbox{is} \; C^{1}, \; u(t) \in T_{c(t)}Q, \; u(0)=0 \; \hbox{and} \;
u(1)=0 \} \, . 
\]

Now, we introduce a special subset $\tilde{\cal C}^{2}(x,y)$ of  ${\cal C}^{2}(x,y)$ which consists of those curves which are in the constraint submanifold $M$
\[
\tilde{\cal C}^{2}(x,y) = \{ {c} \in {\cal C}^{2}(x,y) \, | \;
\dot{c}(t)\in M_{{c}(t)}=M\cap \tau_Q^{-1}(c(t))\, , \; \forall t \in [0,1] \} \, .
\]

Let us consider  the functional ${\cal J}$ defined by
\begin{eqnarray*}
{\cal J} & : & {\cal C}^{2}(x,y)  \longrightarrow \R\\
&& c \mapsto \int_{0}^{1} L(\dot{c}(t)) \, dt \, ,
\end{eqnarray*}
which we want to extremize among the curves satisfying the constraints imposed by $M$, $c \in \tilde{\cal C}^2(x,y)$.

\begin{definition}
A curve ${c} \in \tilde{\cal C}^{2}(x,y)$ will be a {\bf solution of the vakonomic problem} if ${c}$ is a critical point of ${\cal J}_{| \tilde{\cal C}^{2}(x,y)}$.
\end{definition}

Therefore, ${c}$ is a solution of the vakonomic problem if and only if
$d{\cal J}({c})\cdot{u}=0$, for all ${u} \in T_{{c}} \, \tilde{\cal C}^{2}(x,y)$.

\begin{remark}
{\rm In this paper, we will assume that the solution curves $c\in \tilde{\cal C}(x,y)$ admit nontrivial variations in $\tilde{\cal C}(x,y)$. These solutions are called normal in the literature, in opposition to the abnormal ones, which are pathological curves which do not admit nontrivial variations \cite{arnold}. Several investigators have shown the existence of $C^1$, stable under perturbations abnormal solutions \cite{sussmann,montgomery}.
}
\end{remark}

Now, we find a characterization for a curve $c$ to be critical for the vakonomic problem.
\begin{proposition}\label{qwe}
A curve ${c}\in \tilde{\cal C}^{2}(x,y)$ is a normal solution of the vakonomic problem if and only if there exists $\mu: [0,1]\rightarrow \R^m$ such that 
\begin{equation}\label{asdm}
\left\{
\begin{array}{l}
\displaystyle{\frac{d}{dt}\left(\frac{\partial \tilde{L}}{\partial \dot{q}^a}\right)
-\frac{\partial \tilde{L}}{\partial q^a}=
\mu_{\alpha}\left[\frac{d}{dt}\left(\frac{\partial \Psi^{\alpha}}{\partial \dot{q}^a}\right)-\frac{\partial \Psi^{\alpha}}{\partial q^a}\right]+\dot{\mu}_{\alpha}\frac{\partial \Psi^{\alpha}}{\partial \dot{q}^a}},\\
\displaystyle{\dot{\mu}_{\alpha}=
\frac{\partial \tilde{L}}{\partial q^{\alpha}}
-\mu_{\beta}\frac{\partial \Psi^{\beta}}
{\partial q^{\alpha}},}\\
\dot{q}^{\alpha}=\Psi^{\alpha}(q^A, \dot{q}^{a}),
\end{array}
\right.
\end{equation}
where $\tilde{L}: M\rightarrow \R$ is the restriction of $L$ to $M$.
\end{proposition}

\proof{
The condition for a curve to be critical is
\[
0=d{\cal J}(c)\cdot u=\frac{d}{ds}{\cal J}(c_s)\Big|_{s=0}\; ,
\]
for any variation $c_s$ in $\tilde{\cal C}^{2}(x,y)$ of $c$, where $\displaystyle{u=\frac{dc_s}{ds}\Big|_{s=0}}$.

Then, we have that
\begin{eqnarray*}
0=\frac{d}{ds}{\cal J}(c_s)\Big|_{s=0}
&=&\frac{d}{ds}
\left( 
\int^1_0 
L(\dot c_s(t))\, dt
\right)
\Big|_{s=0}\\
&=&\int^1_0 \frac{d}{ds} L(\dot c_s(t))\Big|_{s=0}\, dt\; .
\end{eqnarray*}
In local coordinates, we obtain 
\begin{eqnarray}\label{qqqq}
0&=&
\int^1_0\left(
\frac{\partial L}{\partial q^A}{u}^A+\frac{\partial L}{\partial \dot{q}^a}\dot{u}^a
+\frac{\partial L}{\partial \dot{q}^{\alpha}}
\frac{\partial \Psi^{\alpha}}{\partial q^A}u^A+
\frac{\partial L}{\partial \dot{q}^{\alpha}}
\frac{\partial \Psi^{\alpha}}{\partial \dot{q}^a}\dot{u}^a\right)\, dt\nonumber\\
&=&
\int^1_0\left(\left[
\frac{\partial L}{\partial q^A}
+\frac{\partial L}{\partial \dot{q}^{\alpha}}
\frac{\partial \Psi^{\alpha}}{\partial q^A}
\right]u^A
+
\left[
\frac{\partial L}{\partial \dot{q}^a}
+
\frac{\partial L}{\partial \dot{q}^{\alpha}}
\frac{\partial \Psi^{\alpha}}{\partial \dot{q}^a}
\right]\dot{u}^a
\right)\, dt\\
&=&
\int^1_0\left(
\frac{\partial \tilde{L}}{\partial q^A}u^A
+
\frac{\partial \tilde{L}}{\partial \dot{q}^a}\dot{u}^a
\right)\, dt \nonumber
\end{eqnarray}
>From (\ref{constraints}) we know that the infinitesimal variations $u^A$, $1\leq A\leq n$, are not arbitrary. Consider the functions $\mu_{\alpha}$ defined as the solutions of the following system of first order differential equations
\[
\dot{\mu}_{\alpha}=
\frac{\partial \tilde{L}}{\partial q^{\alpha}}\Big|_{c}
-\mu_{\beta}\frac{\partial \Psi^{\beta}}{\partial q^{\alpha}}\Big|_{c} \, , \; 1 \le \alpha \le m \, .
\]
Then, using the fact that $\displaystyle{
\dot{u}^{\alpha}= \frac{\partial \Psi^{\alpha}}{\partial q^A}u^A+\frac{\partial \Psi^{\alpha}}{\partial \dot{q}^a}\dot{u}^a}$, we get 
\begin{eqnarray*}
\frac{d}{dt}(\mu_{\alpha}u^{\alpha})&=&
\mu_{\alpha}\dot{u}^{\alpha}+\left(
\frac{\partial \tilde{L}}{\partial q^{\alpha}}
-\mu_{\beta}\frac{\partial \Psi^{\beta}}{\partial q^{\alpha}}\right)u^{\alpha}\\
&=&
u^{\alpha}
\frac{\partial \tilde{L}}{\partial q^{\alpha}}
+\mu_{\alpha}\frac{\partial \Psi^{\alpha}}{\partial q^a}
u^a+\mu_{\alpha}\frac{\partial \Psi^{\alpha}}{\partial \dot q^a}\dot{u}^a\; ,\\
\end{eqnarray*}
or, equivalently, 
\[
u^{\alpha}
\frac{\partial \tilde{L}}{\partial q^{\alpha}}
=\frac{d}{dt}(\mu_{\alpha}u^{\alpha})-
\mu_{\alpha}\frac{\partial \Psi^{\alpha}}{\partial q^a}
u^a-\mu_{\alpha}\frac{\partial \Psi^{\alpha}}{\partial \dot q^a}\dot{u}^a\; .
\]
Substituting the last expression in (\ref{qqqq}) and integrating by parts, we obtain
\[
d{\cal J}(c)\cdot u=
\int^1_0
\left[
\frac{\partial \tilde{L}}{\partial q^a}
-
\mu_{\alpha}\frac{\partial \Psi^{\alpha}}{\partial q^a}
\right]u^a
+
\left[
\frac{\partial \tilde{L}}{\partial \dot{q}^a}
-\mu_{\alpha}\frac{\partial \Psi^{\alpha}}{\partial \dot q^a}
\right]\dot{u}^a
\, dt \, .
\]
Now, since
\begin{eqnarray*}
\left[
\frac{\partial \tilde{L}}{\partial \dot{q}^a}
-\mu_{\alpha}\frac{\partial \Psi^{\alpha}}{\partial \dot q^a}
\right]\dot{u}^a
=
\frac{d}{dt}\left(\left[
\frac{\partial \tilde{L}}{\partial \dot{q}^a}
-\mu_{\alpha}\frac{\partial \Psi^{\alpha}}{\partial \dot q^a}\right]u^a\right)
-\frac{d}{dt}\left(
\frac{\partial \tilde{L}}{\partial \dot{q}^a}
-\mu_{\alpha}\frac{\partial \Psi^{\alpha}}{\partial \dot q^a}
\right)u^a\; ,
\end{eqnarray*}
using again integration by parts, we can write
\[
0= \int^1_0 
\left[
\frac{\partial \tilde{L}}{\partial q^a}
-
\mu^{\alpha}\frac{\partial \Psi^{\alpha}}{\partial {q}^a}
-\frac{d}{dt}\left(
\frac{\partial \tilde{L}}{\partial \dot{q}^a}
-\mu^{\alpha}\frac{\partial \Psi^{\alpha}}{\partial \dot{q}^a}
\right)
\right]u^a
\, dt\; .\\
\]
As the infinitesimal variations $u^a$ are arbitrary, the fundamental lemma of the Calculus of Variations applies and we can assert that $d{\cal J}(c)\cdot u=0$ if and only if $c$ and $\mu_{\alpha}$ satisfy equations (\ref{asdm}).} \QED

\begin{remark}\label{exten}
{\rm 
The usual way in which the equations of motion for vakonomic mechanics are presented is the following
\begin{equation}\label{asdo}
\left\{
\begin{array}{l}
\displaystyle{ \frac{d}{dt}\left( \frac{\partial L}{\partial \dot{q}^A} \right) -
\frac{\partial L}{\partial q^A}  = \dot{\lambda}_{\alpha} \frac{\partial \Phi^{\alpha}}{\partial \dot q^A}
 + \lambda_{\alpha} \left[ \frac{d}{dt}\left( \frac{\partial \Phi^{\alpha}}{\partial \dot{q}^A} \right) -
\frac{\partial \Phi^{\alpha}}{\partial q^A}\right]}\,,\\
\\
\Phi^{\alpha}(q, \dot{q}) = 0, \ 1 \leq \alpha \leq m \,,
\end{array}
\right.
\end{equation}

where $\Phi^{\alpha}=\Psi^{\alpha}-\dot{q}^{\alpha}$ and 
$\displaystyle{
\lambda_{\alpha}=\frac{\partial L}{\partial \dot q^{\alpha}}-\mu_{\alpha},\ 1\leq \alpha \leq m\; .}$
Observe that, in contrast to equations (\ref{asdm}), equations (\ref{asdo}) are expressed in terms of the ambient lagrangian $L: TQ\rightarrow \R$. Equations (\ref{asdm}) stress how the information given by $L$ outside $M$ is irrelevant to obtain the vakonomic equations, contrary to what happens in nonholonomic mechanics (see  Section \ref{pppp} below).

Equations (\ref{asdo}) can be seen as the Euler-Lagrange equations for the extended lagrangian ${\cal L}=L+\lambda_{\alpha}\Phi^{\alpha}$. We will not follow this approach here, which has been exploited fruitfully in \cite{Fa,KoMa,MaCoLe}. Finally, note that if we consider the extended lagrangian $\lambda_0 L+\lambda_{\alpha}\Phi^{\alpha}$, with $\lambda_0=0$ or $1$, then we recover all the solutions, both the normal and the abnormal ones \cite{arnold}.
}
\end{remark} 

\section{Geometric approach to vakonomic mechanics}\label{S2}

We will develop a geometric characterization of vakonomic mechanics following an approach similar to the formulation given by Skinner and Rusk \cite{SR,SR2} for singular lagrangians (see also \cite{clr,IbMa,LR,CL}). This characterization is specially interesting, since enables us to study both linear and nonlinear constraints in an intrinsic way. We will show in the examples the utility of this formulation.

Consider the Whitney sum of $T^*Q$ and $TQ$, $T^*Q\oplus TQ$,
and its canonical projections
\[
pr_1: T^*Q \oplus TQ\longrightarrow T^*Q \, , \quad pr_2: T^*Q \oplus TQ\longrightarrow TQ \, .
\]
Let us take the submanifold $W_0=pr_2^{-1}(M)$, where $M$ is the constraint submanifold, locally determined by the constraint equations $\Phi^{\alpha}=0$, $1\leq \alpha \leq m$. We will denote $W_0=T^*Q\times_Q M$ and $\pi_1={pr_1}_{\big|{W_0}}$, $\pi_2={pr_2}_{\big|{W_0}}$.
Now, define on $T^*Q\times_Q M$ the presymplectic 2-form $\omega=\pi_1^*\omega_Q$, where $\omega_Q$ is the canonical symplectic form on $T^*Q$. Observe that the rank of this presymplectic form is equal to $2n$ everywhere. Define also the function
\[
H_{W_0}=\langle \pi_1,\pi_2\rangle -\pi_2^* \tilde{L}\; ,
\]
where $\tilde{L}: M\rightarrow \R$ is the restriction of $L$ to $M$.

If $(q^A)$ are local coordinates on a neighborhood $U$ of $Q$, $(q^A, \dot{q}^a)$ coordinates on $TU\cap M$ and $(q^A, p_A)$ the induced coordinates on $T^*U$, then we have induced coordinates
$(q^A, p_A, \dot{q}^a)$ on $T^*U\times_Q (TU\cap M)$. Locally, the hamiltonian function $H_{W_0}$ reads as
\[
H_{W_0}(q^A, p_A, \dot{q}^a)=p_a\dot{q}^a+p_{\alpha}\Psi^{\alpha}-\tilde{L}(q^A, \dot{q}^a) \, ,
\]
and the 2-form $\omega$ is $\omega=dq^A\wedge dp_A$.

Now, we will see how the dynamics of the vakonomic system (\ref{asdm}) is determined by studying the solutions of the equation 
\begin{equation}\label{vhs}
i_X\omega=dH_{W_0} \, .
\end{equation}
Thus, we are justified to employ the following terminology:
\begin{definition}
The presymplectic hamiltonian system  $(T^*Q\times_Q M,\omega, H_{W_0})$ will be called {\bf  vakonomic hamiltonian system.}
\end{definition}
Being the system $(T^*Q\times_Q M,\omega, H_{W_0})$ presymplectic, we apply to it  the Gotay-Nester's constraint algorithm \cite{Got,GN}. First we consider the points $W_1$ of $T^*Q\times_Q M$ where  (\ref{vhs}) has a solution. This first constraint submanifold is determined by
\[
W_1=\{ x\in T^*Q\times_Q M\; | \; dH_{W_0}(x)(V)=0,\ \forall V\in \ker \omega(x)\}\; .
\]
Locally, $\ker \omega=\displaystyle{\langle \frac{\partial}{\partial \dot{q}^a}\rangle} $. Therefore, the constraint submanifold $W_1$ is locally characterized  by the vanishing of the constraints
\[
\varphi_a=p_a+p_{\alpha}\frac{\partial \Psi^{\alpha}}{\partial \dot{q}^a}-
\frac{\partial \tilde{L}}{\partial \dot{q}^a}=0\; ,
\]
or, equivalently,
\begin{equation}\label{asz}
p_a=\frac{\partial \tilde{L}}{\partial \dot{q}^a}
-p_{\alpha}\frac{\partial \Psi^{\alpha}}{\partial \dot{q}^a}\, , \; m+1 \le a \le n \, .
\end{equation}

Expanding the expressions in equation (\ref{vhs}) and equating coefficients, we obtain that the equations of motion along $W_1$ are
\begin{eqnarray*}
\dot{q}^{A}&=&\frac{\partial H_{W_0}}{\partial p_A}\; ,\\
\dot{p}_A&=&-\frac{\partial H_{W_0}}{\partial q^A}\; ,
\end{eqnarray*}
which is equivalent to
\begin{equation}\label{asx}
\left\{
\begin{array}{rcl}
\dot{q}^{\alpha}&=&\Psi^{\alpha}(q^A, \dot{q}^a)\; ,\\
\displaystyle{\dot{p}_{\alpha}}&=& \displaystyle{\frac{\partial \tilde{L}}{\partial q^{\alpha}}-p_{\beta}\frac{\partial \Psi^{\beta}}{\partial q^{\alpha}}}\; ,\\
\displaystyle{\frac{d}{dt}
\left(\frac{\partial \tilde{L}}{\partial \dot{q}^a}
-p_{\alpha}\frac{\partial \Psi^{\alpha}}{\partial \dot{q}^a}
\right)}
&=&\displaystyle{\frac{\partial \tilde{L}}{\partial q^{a}}-p_{\beta}\frac{\partial \Psi^{\beta}}{\partial q^{a}}} \; .
\end{array}
\right.
\end{equation}
Observe that these equations are precisely the vakonomic equations of motion (\ref{asdm}), where now $p_{\alpha}=\mu_{\alpha}$, $1\leq \alpha\leq m$.  

\begin{remark}
{\rm 
The momenta $p_{\alpha}$, $1\leq \alpha\leq m$, play the role of the Lagrange multipliers, but they do not have any physical meaning (see \cite{taimanov}).
}
\end{remark}

Therefore, a vector field $X$ solution of equation (\ref{vhs}), will be of the form
\begin{eqnarray*}
X&=&\dot{q}^a \left( \frac{\partial }{\partial q^a} + \left( \frac{\partial^2 \tilde{L}}{\partial q^a \partial \dot{q}^b} - p_\gamma \frac{\partial^2 \Psi^\gamma}{\partial q^a \partial \dot{q}^b} \right) \frac{\partial}{\partial p_b} \right) + \Psi^\alpha \left( \frac{\partial }{\partial q^\alpha} + \left( \frac{\partial^2 \tilde{L}}{\partial q^\alpha \partial \dot{q}^b} - p_\gamma \frac{\partial^2 \Psi^\gamma}{\partial q^\alpha \partial \dot{q}^b} \right) \frac{\partial}{\partial p_b} \right) \\
&& + \, \bar{Y}^a \left( \frac{\partial}{\partial \dot{q}^a} + \left( \frac{\partial^2 \tilde{L}}{\partial \dot{q}^a \partial \dot{q}^b} - p_\gamma \frac{\partial^2 \Psi^\gamma}{\partial \dot{q}^a \partial \dot{q}^b} \right) \frac{\partial }{\partial p_b} \right) + \left( \frac{\partial \tilde{L}}{\partial q^{\alpha}}-p_{\beta}\frac{\partial \Psi^{\beta}}{\partial q^{\alpha}}
\right) \left( \frac{\partial}{\partial p_\alpha} - \frac{\partial \Psi^\alpha}{\partial \dot{q}^b}\frac{\partial}{\partial p_b} \right) \, .
\end{eqnarray*}

\begin{figure}[ht]
\centering
\fbox{\begin{minipage}{10cm}
\vspace{0.4cm}
\centering
\setlength{\unitlength}{1cm}
\begin{center}
\begin{picture}(6.2,6.2)(-1,0)
\put(-2,3){\makebox(0,0)[r]{$M$}}
\put(6,3){\makebox(0,0)[l]{$T^*Q$}}
\put(2,0){\makebox(0,0)[c]{$Q$}}
\put(2,6){\makebox(0,0)[c]{$T^*Q\oplus TQ$}}
\put(2,4){\makebox(0,0)[c]{$W_0=T^*Q\times_Q M$}}
\put(2,2){\makebox(0,0)[c]{$W_1$}}
\put(2,4.3){\vector(0,1){1.4}}
\put(2,2.3){\vector(0,1){1.4}}
\put(2,1.7){\vector(0,-1){1.4}}
\put(0,4){\vector(-3,-1){2}}
\put(4,4){\vector(3,-1){2}}
\put(-1,4){\makebox(0,0)[r]{$\pi_2$}}
\put(5,4){\makebox(0,0)[l]{$\pi_1$}}
\put(-2,2.6){\vector(3,-2){3.5}}
\put(6,2.6){\vector(-3,-2){3.5}}
\put(-0.6,1.5){\makebox(0,0)[r]{$(\tau_Q)_{|M}$}}
\put(3.4,1.5){\makebox(0,0)[l]{$\pi_Q$}}
\end{picture}
\end{center}
\vskip1cm
\end{minipage}}
\centering
\vspace{0.2cm}

Geometric formulation of vakonomic mechanics
\end{figure}

Nevertheless, the solutions on $W_1$ may not be tangent to $W_1$. In such a case, we have to restrict $W_1$ to the submanifold $W_2$ where these solutions are tangent to $W_1$. Proceeding further, we obtain a sequence of submanifolds (we are assuming that all the subsets obtained are submanifolds)
\[
\cdots\hookrightarrow  W_k\hookrightarrow  \cdots \hookrightarrow  W_2\hookrightarrow  W_1\hookrightarrow  W_0=T^*Q\times_Q M\; .
\]
Algebraically, these constraint submanifolds may be described as
\[
W_i=\{x\in T^*Q\times_Q M\; | \;
dH_{W_0}(x)(v)=0\; , \  \forall v\in T_x W_{i-1}^{\perp} \;  \}\;  , \quad i\geq 1\; ,
\]
where
\[
T_x W_{i-1}^{\perp} = \{ v \in T_x(T^*Q\times_Q M) \; | \;  \omega(x)(u,v) = 0 \; ,\
\forall u \in T_x W_{i-1} \; \} \; .
\]
If this constraint algorithm stabilizes, i.e., if there exists a positive integer $k \in \N$ such that $W_{k+1}=W_k$ and $\dim W_k\not=0$, then we will have obtained a final constraint submanifold $W_f=W_k$ on which a vector field $X$ exists such that
\[
\left(i_X\omega=dH_{W_0}\right)_{| W_f}\; .
\]

Note that on $W_f$ we will have an explicit solution of the vakonomic dynamics. A very important particular case is when the final constraint submanifold is the first one, $W_f=W_1$. Observe that the dimension of $W_1$ is even, $\dim W_1=2n$. In the sequel, we will investigate when this constraint submanifold is equipped with a symplectic 2-form in order to determine an unique solution of the vakonomic equations.
Obviously, this geometrical study is related with the explicit or implicit character of the second order differential equations obtained in (\ref{asdm}).

Denote by $\omega_{W_1}$ the restriction of the presymplectic $2$-form $\omega$ to $W_1$.

\begin{proposition}\label{axc}
$(W_1, \omega_{W_1})$ is a symplectic manifold if and only if, for any choice of coordinates $(q^A,p_A, \dot{q}^a)$ on $T^*Q\times_Q M$, 
\[
\det\left(\frac{\partial^2 \tilde{L}}{\partial \dot{q}^a\partial \dot{q}^b}
-p_{\alpha}\frac{\partial^2 \Psi^{\alpha}}{\partial \dot{q}^a \partial \dot{q}^b}\right)\not=0
\]
for all point in $W_1$.
\end{proposition}
{\bf Proof:} $\omega_{W_1}$ is symplectic if and only if 
\begin{equation}\label{aw}
T_xW_1\cap (T_xW_1)^{\perp}=0 \, ,
\end{equation}
for all $x\in W_1$.
Condition (\ref{aw}) is satisfied if and only if the matrix
$\displaystyle{
d\varphi_a(\frac{\partial}{\partial p_b})
}$
is regular, that is, 
\[
\det\left(\frac{\partial^2 \tilde{L}}{\partial \dot{q}^a\partial \dot{q}^b}
-p_{\alpha}\frac{\partial^2 \Psi^{\alpha}}{\partial \dot{q}^a \partial \dot{q}^b}\right)\not=0\; ,
\]
for all $x\in W_1$.
\QED

In this case, the equations of motion (\ref{asx})
are rewritten as the following system of algebraic and explicit differential equations

\parbox{15cm}{
\[
\left\{
\begin{array}{l}
\dot{q}^{\alpha}=\Psi^{\alpha}(q^A, \dot{q}^a)\; ,\\
\displaystyle{\dot{p}_{\alpha}= \frac{\partial \tilde{L}}{\partial q^{\alpha}}-p_{\beta}\frac{\partial \Psi^{\beta}}{\partial q^{\alpha}}}\; ,\\
\ddot{q}^a=
\displaystyle{-\bar{\cal {C}}^{ab}\left[
\dot{q}^A\frac{\partial^2 \tilde{L}}{\partial q^A\partial\dot{q}^b}-
\dot{q}^A p_{\alpha}
\frac{\partial^2 \Psi^{\alpha}}{\partial q^A \partial \dot{q}^b}\right.}\\
\qquad \displaystyle{\left. -\frac{\partial \tilde{L}}{\partial q^b}
+p_{\alpha}\frac{\partial \Psi^{\alpha}}{\partial q^b}
-
\left(
\frac{\partial \tilde{L}}{\partial q^{\gamma}}-p_{\beta}\frac{\partial \Psi^{\beta}}{\partial q^{\gamma}}\right)
\frac{\partial \Psi^{\gamma}}{\partial \dot{q}^b}\right]}\; ,
\end{array}
\right.
\]
}\
\parbox{1.45cm}{
\begin{equation}\label{ytr}
\end{equation}}
where
\begin{equation}\label{5}
\bar{\cal C}_{ab}=\frac{\partial^2 \tilde{L}}{\partial \dot{q}^a\partial \dot{q}^b}-p_{\alpha}
\frac{\partial^2 \Psi^{\alpha}}{\partial \dot{q}^a \partial \dot{q}^b} \, ,
\end{equation}
and $(\bar{\cal C}^{ab})$ denotes the inverse matrix of $(\bar{\cal C}_{ab})$.

\begin{remark}
{\rm The characterization found in Proposition \ref{axc} for the symplecticness of the manifold $(W_1,\omega_{W_1})$ implies that the constraint equations
\[
\varphi_a=p_a+p_{\alpha}\frac{\partial \Psi^{\alpha}}{\partial \dot{q}^a}-
\frac{\partial \tilde{L}}{\partial \dot{q}^a}=0 \; , m+1 \le a \le n \, ,
\]
define locally the variables $\dot{q}^a$, $m+1 \le a \le n$, by the implicit function theorem. That is, we have
\[
\dot{q}^a = \varsigma^a(q^A,p_A) \, , \; m+1 \le a \le n \, .
\]
Therefore, we can also consider local coordinates $(q^A,p_A)$ on $W_1$. In such a case, the symplectic form has the following local expression
\[
\omega_{W_1} = dq^A \wedge dp_A \, ,
\]
and the restriction of the hamiltonian $H_{W_0}$ to $W_1$
\[
H_{W_1} = p_a \varsigma^a + p_\alpha \Psi^\alpha - \bar{L}(q^A,p_A) \, ,
\]
where $\bar{L}(q^A,p_A)= \tilde{L}(q^A,\varsigma^a(q^A,p_A))$. Consequently, equations (\ref{ytr}) are rewritten in hamiltonian form as
\begin{equation}
\left\{
\begin{array}{rcl}
\dot{q}^A &=& \displaystyle{\frac{\partial H_{W_1}}{\partial p_A}} \\
\dot{p}_A &=& -\displaystyle{\frac{\partial H_{W_1}}{\partial q^A}}
\end{array}
\right.
\end{equation}
This choice of coordinates is very common in optimal control theory.
}
\end{remark}

Now, observe that, if the constraints are linear on the velocities, 
we can write
\[
\dot{q}^{\alpha}=\Psi^{\alpha}_a(q)\dot{q}^a\; .
\]
Then, from Proposition \ref{axc}, 
$\omega_{W_1}$ is symplectic if and only if
\[
\det\left(\frac{\partial^2 \tilde{L}}{\partial \dot{q}^a\partial \dot{q}^b}
\right)\not=0\; .
\]

\begin{proposition}\label{asdd}
Suppose that the constraints are given by
\[
\dot{q}^{\alpha}=\Psi^{\alpha}_a(q)\dot{q}^a\, ,\ 1\leq \alpha\leq m,
\]
and the lagrangian $L$ is regular. 
Denote by $(W^{AB})$ the inverse matrix of the Hessian matrix of $L$.
In this case, 
$\omega_{W_1}$ is symplectic on $W_1$ if and only if 
the constraints are compatible, that is,
the matrix whose entries are 
\[
{\cal C}^{\alpha \beta}= W^{ab}\Psi^{\alpha}_a\Psi^{\beta}_b
-W^{\alpha b}\Psi^{\beta}_b
-W^{a \beta}\Psi^{\alpha}_a+W^{\alpha\beta}\, ,
\]
is nonsingular. 
\end{proposition}
\proof{
See the geometrical proof of Theorem IV. 3 in reference \cite{MD}.\QED
}

\begin{remark}
{\rm 
The compatibility condition guarantees the existence and unicity of the solutions for the nonholonomic problem with lagrangian $L$ and constraint submanifold $M$ \cite{MD,WSF}. 
}
\end{remark}

\subsection{Applications to economy}

The variational calculus is an indispensable tool in many classical and recent economic papers \cite{KaHa,Mag,Ra,Sa,SaRa}. In fact, a typical optimization problem in modern economics deals with the problem of maximizing or minimizing the functional 
\[
\int^T_0 D(t)U[f(t,k,\dot{k})]\; dt 
\]
subject or not to constraints.
Here, $D(t)$ is a discount rate factor, $U$ an utility function, $f$ a consumption function and  $k$ the capital labor-ratio.
It is usual to find dynamical economic models with nonholonomic constraints. For instance, the revision of the expected rate of inflaction may be expressed in terms of the nonholonomic constraint
\[
\dot{\pi}=j(p-\pi), \ 0<j\leq 1
\]
where $\pi$ and $p$ are the expected and actual rates of inflaction, respectively.

In economics, it is also very common to deal with an explicit dependence of the time. We plan to extend the geometric formulation of vakonomic dynamics to the non-autonomous case in a forthcoming paper.

\begin{example}[Closed von Neumann System \cite{Sa,SaRa}]
{\rm 
Consider the transformation function which relates $n$ capital goods $K_1, K_2,\ldots, K_n$ and the net capital formations $\dot{K}_1,\dot{K}_2,\ldots, \dot{K}_n$
\[
F(K_1,\ldots, K_n, \dot{K}_1,\ldots, \dot{K}_n)=
K_1^{\alpha_1}K_2^{\alpha_2}\cdots K_n^{\alpha_n} - \left[ \dot{K}_1^{2}+\ldots+\dot{K}_n^{2}\right]^{1/2}\; ,
\]
with $\alpha_1+\alpha_2+\ldots +\alpha_n=1$.
The von Neumann problem is to maximize
\[
\int^T_0 \dot{K}_n\; dt\quad \hbox{ subject to }\quad F(K_1,\ldots, K_n, \dot{K}_1,\ldots, \dot{K}_n)=0
\]
and appropiate initial conditions.

Applying our formulation it is possible to write this problem as a presymplectic system on $W_0=R^{3n-1}$. The constraint $F=0$ can be rewritten as
\[
\dot{K}_1 = \left(K_1^{2\alpha_1},\ldots, K_n^{2\alpha_n} - \sum_{i=2}^{n}\dot{K}_i^2\right)^{1/2} = \Psi(K_1,\ldots, K_n,\dot{K}_2,\ldots, \dot{K}_n) \, .
\]
Taking coordinates $(K_1,\ldots, K_n, \dot{K}_2,\ldots, \dot{K}_{n}, P^1,\ldots, P^n)$ we have that
\begin{eqnarray*}
\omega&=&\sum_{i=1}^n dK_i\wedge dP^i\\
H_{W_0}&=&\sum_{i=2}^n P^i\dot{K}_i+P^1
\cdot\left(
K_1^{2\alpha_1} K_2^{2\alpha_2}\cdots K_n^{2\alpha_n}
-\sum_{i=2}^{n}\dot{K}_i^2 \right)^{1/2}-\dot{K}_n
\end{eqnarray*}
Applying the Gotay and Nester algorithm new constraints arise
\begin{eqnarray*}
P^i&=&P^1\dot{K}_i\left(
K_1^{2\alpha_1}K_2^{2\alpha_2}\cdots K_n^{2\alpha_n}
-\sum_{i=2}^{n}\dot{K}_i^2\right)^{-1/2}\; ,\quad 2\leq i\leq n-1\\
P^n&=&1+P^1\dot{K}_n\left(
K_1^{2\alpha_1}K_2^{2\alpha_2}\cdots K_n^{2\alpha_n}
-\sum_{i=2}^{n}\dot{K}_i^2\right)^{-1/2}
\end{eqnarray*}
Therefore, from (\ref{asx}) the initial system is determined solving the following $n$ differential equations on the variables $(K_1,\ldots, K_n, \dot{K}_2, \ldots,\dot{K}_n, P^1)$
\begin{equation}
\left\{
\begin{array}{rcl}
\dot{P}^1 &=&-P^1\alpha_1\left(K_1^{2\alpha_1-1}K_2^{2\alpha_2}\cdots K_n^{2\alpha_n}\right)\cdot
G
\\
0 & = & \dot{P}^1\dot{K}_i G\\
&&+P^1\left[\left(\ddot{K}_i+\alpha_i\left(K_1^{2\alpha_1}
\cdots K_i^{2\alpha_i-1}\cdots K_n^{2\alpha_n}\right)\right)G
+\dot{K}_i\frac{d}{dt}(G)\right]\; ,\ 2\leq i\leq n
\end{array}
\right.
\end{equation}
where 
\[
G(K_1,\ldots, K_n, \dot{K}_2,\ldots, \dot{K}_n)= \frac{1}{\Psi(K_1,\ldots, K_n, \dot{K}_2,\ldots, \dot{K}_n)} \, .
\]
}
\end{example}

\subsection{Principal kinematic locomotion systems}

This kind of systems includes the motion of inchworms, paramecia, mobile vehicles, robotic snakes, etc \cite{KeMu,KeMu2,OsBu2}. The study of the motion relies on the simple fact that the process of locomotion can be divided into internal (shape) variables and position (group) variables. The internal variables are assumed to be directly controlled. The motion of these variables couple to produce a net change in the position and orientation of the moving body. Moreover, locomotion systems are characterized by the fact that the constraints are usually invariant with respect to the group action. Consequently, one is provided with a useful mathematical structure to work with: a principal connection $\gamma$ on a principal fibre bundle, $\pi: Q \longrightarrow Q/G$. The constraint submanifold is precisely the horizontal distribution ${\cal H}$ of $\gamma$.

In the following, we will investigate optimal control problems for these systems. Take local coordinates in $Q$, $(r,g)$, where $r$ are coordinates in the base manifold $Q/G$ and $g$ in the Lie group $G$. Let us assume that we are given a quadratic cost function $C$, locally expressed as $C(r,\dot{r}) = \displaystyle{\frac{1}{2}}C_{a b}(r)\dot{r}^a\dot{r}^b$ and define the functional
\[
{\cal J} = \int^1_0 C(r,\dot{r}) dt \, .
\]
The cost function $C$ depends only on the shape variables, which corresponds to calculating the cost of the control effort. Then, the optimal control problem is to obtain the inputs that will minimize ${\cal J}$, while steering the state from $(r_0,g_0)$ to $(r_1,g_1)$.

>From the vakonomic point of view, we have the lagrangian $L=C$, the constraint submanifold $M={\cal H}$ and the infinite-dimensional manifold $\tilde{\cal C}^2((r_0,g_0),(r_1,g_1))$. We will not go into depth here with the mathematical structure associated to the connection $\gamma$ and the symmetries $G$. This will be a subject of future research.  

\begin{example}[Locomotion at low Reynolds' number \cite{KeMu2,Os2}]
{\rm The kinematic connection for the paramecia can be determined by examining the Stokesian flow around a deformable cylindrical body \cite{KeMu2}. We parametrize the body in polar coordinates $(R,\theta)$ as
\[
R=1 + \epsilon (k_1(t) \cos 2\theta + k_2(t) \cos 3\theta) \, .
\]
The shape variables for this system are $r=(k_1,k_2)$. Let $x$ denote the motion of the centroid of the body in the direction given by $\theta=0$. Symmetry arguments show that all resultant motion must be directed along this ray. Therefore, we have a principal bundle $\pi: \R^3 \longrightarrow \R^2$ with structure group $G=\R$. In \cite{KeMu2}, it is shown that the viscous connection can be written to first order as
\[
\gamma = \dot{x} + \frac{\epsilon^2}{4}(k_2\dot{k}_1+2k_1 \dot{k}_2) \, .
\]
We consider the optimal control problem associated with the simple quadratic cost function $C=\dot{k}_1^2 + \dot{k}_2^2$.

Taking coordinates $(k_1,k_2,x,\dot{k_1},\dot{k_2},p_1,p_2,p_3)$ in $W_0=M \times_Q T^*Q$, we have that
\begin{eqnarray*}
\omega &=& dk_1 \wedge dp_1 + dk_2 \wedge dp_2 + dx \wedge dp_3 \\
H_{W_0} &=& p_1 \dot{k}_1 + p_2 \dot{k}_2 - p_3 \frac{\epsilon^2}{4}(k_2\dot{k}_1+2k_1 \dot{k}_2) - \dot{k}_1^2 - \dot{k}_2^2  \, .
\end{eqnarray*}
Applying the Gotay-Nester algorithm, new constraints arise
\begin{eqnarray*}
p_1 &=& \frac{\epsilon^2}{4}k_2 p_3 + 2 \dot{k}_1 \, , \\
p_2 &=& \frac{\epsilon^2}{2}k_1 p_3 + 2 \dot{k}_2 \, .  
\end{eqnarray*} 
Now, equations (\ref{asx}) read as
\begin{equation}
\left\{
\begin{array}{rcl}
\dot{x}&=&-\displaystyle{\frac{\epsilon^2}{4}(k_2\dot{k}_1+2k_1 \dot{k}_2)} \\
\dot{p}_3 &=& 0 \\
\ddot{k}_1 &=& p_3 \displaystyle{\frac{\epsilon^2}{8}} \dot{k}_2 \\
\ddot{k}_2 &=& -p_3 \displaystyle{\frac{\epsilon^2}{8}} \dot{k}_1
\end{array}
\right.
\end{equation}
These equations can be easily integrated. Setting $\displaystyle{a=p_3 \frac{\epsilon^2}{8}}$, we have
\begin{eqnarray*}
k_1(t) &=& B \cos at + C \sin at - \frac{A}{a} \\
k_2(t) &=& B \sin at - C \cos at + D 
\end{eqnarray*}
for some constants $A$, $B$, $C$, $D$, which is the same result obtained in \cite{Os2}.
}
\end{example}

\subsection{Vakonomic mechanics and sub-riemannian geometry
}
Let $Q$ be an $n$-dimensional manifold with a smooth distribution ${\cal D}$ of constant rank $n-m$. A {\bf sub-riemannian metric} \cite{Br} on ${\cal D}$ is a smoothly varying in $q$ positive definite quadratic form $g_q$ on ${\cal D}_q$. A piecewise smooth curve $\gamma$ in $Q$ is called admissible if $\dot{\gamma}(t) \in {\cal D}$. We define the length of such a curve in the usual way
\[
l(\gamma)=\int \sqrt{g(\dot{\gamma}(t),\dot{\gamma}(t))}\, dt\; .
\]

The sub-riemannian distance between two points $x, y\in Q$ is defined as
\[
d(x,y)=\hbox{inf}_{\gamma}(l(\gamma)) \; ,
\]
for all admissible curves $\gamma$ connecting $x$ and $y$. The distance is taken infinite if there is no such a path.

A curve which realizes the distance between two points is called a {\bf minimizing geodesic}. It is easy to show that $\gamma$ is a minimizing geodesic if it minimizes the functional
\[
\int g(\dot{\gamma}(t), \dot{\gamma}(t))\,dt\; ,
\]
among all the admissible curves with the same endpoints.

Let $\mu_1,\ldots, \mu_m$ be a basis of $1$-forms for the annihilator ${\cal D}^o$. Then an admissible path must verify the nonholonomic constraints
\begin{equation}\label{sub-1}
\mu_i (\dot{\gamma})=0,\ 1\leq i\leq m \; .
\end{equation}
Thus, we see that the problem of finding minimizing geodesics corresponds exactly to the problem of solving the vakonomic problem determined by the restricted lagrangian $\tilde{L}=\frac{1}{2}g$ and the nonholonomic constraints (\ref{sub-1}). Note that, as the constraints are linear and $g$ is positive definite, the Gotay-Nester's algorithm always ends in the first step, $W_f=W_1$.

\begin{example}[Sub-riemannian geometry: The Martinet case \cite{BoCh}]
{\rm
Let $U$ be an open set of $\R^3$ containing $0$ and ${\cal D}$ the distribution on $U$ determined by annihilation of the Martinet $1$-form
\[
\mu=dz-\frac{y^2}{2}\,dx\; ,
\]
which determines the constraint function $\displaystyle{\dot{z}=\frac{y^2}{2}\dot{x}}$. Moreover, consider the restricted lagrangian
\[
\tilde{L}(x,y,z,\dot{x},\dot{y})=\frac{1}{2}(\dot{x}^2+\dot{y}^2)\; .
\]
This corresponds to the flat metric case in \cite{BoCh}.

The normal minimizing geodesics for this problem are determined by solving the presymplectic hamiltonian system 
\[
i_X\omega=dH_{W_0},
\]
where, locally,
\begin{eqnarray*}
\omega&=&dx\wedge dp_x+dy\wedge dp_y+dz\wedge dp_z\; ,\\
H_{W_0}&=& \dot{x}p_x+\dot{y}p_y+\frac{y^2}{2}\dot{x}p_z-\frac{1}{2}(\dot{x}^2+\dot{y}^2)\; .
\end{eqnarray*}
Applying the Gotay-Nester constraint algorithm we obtain the new constraints
\begin{eqnarray*}
p_x&=&\dot{x}-\frac{y^2}{2}p_z\; ,\\
p_y&=&\dot{y}
\end{eqnarray*}
In this particular case, the equations of motion (\ref{asx}), in coordinates $(x,y,z, \dot{x},\dot{y},p_z)$, are
\begin{equation}
\left\{
\begin{array}{rcl}
\dot{z}&=&\displaystyle{\frac{y^2}{2}\dot{x}}\\
\dot{p}_z&=&0\\
\displaystyle{\frac{d}{dt}\left(\dot{x}-\frac{y^2}{2}p_z\right)}&=&0\\
\ddot{y}&=&-y\dot{x}p_z
\end{array}
\right.
\end{equation}
which are obviously integrable by quadratures (compare with \cite{BoCh}).
}
\end{example} 

\section{Vershik-Gershkovich and vakonomic hamiltonian approaches compared} 

In the precedent section, we have found an intrinsic geometric approach to vakonomic dynamics. It is possible to give an alternative geometric formulation of the vakonomic equations of motion  related to the one of Vershik and Gershkovich \cite{VG}. A key element to obtain this alternative description will be the next fibred diffeomorphism
\[
\begin{array}{rrcl}
F:& T^*Q\oplus TQ&\longrightarrow& T^*Q\oplus TQ\\
  & (\alpha, v)&\longmapsto      & (\alpha-Leg_L(v), v) \, ,
\end{array}
\]
for any $\alpha\in T_x^*Q$, $v\in T_xQ$ and $x\in Q$. Here, $Leg_L: TQ\rightarrow T^*Q$ denotes the Legendre transformation associated to the lagrangian $L$, which in local coordinates reads as $Leg_L(q^A,\dot{q}^A)=(q^A,\displaystyle{\frac{\partial L}{\partial \dot{q}^A}})$. It is clear that $F(T^*Q\times_Q M)=T^*Q\times_Q M$. We will see how in the case of linear constraints, we ``recover" the Vershik-Gershkovich formulation. As a by-product, we will have obtained a generalization of their formulation to the case of nonlinear constraints.

Consider on $T^*Q\oplus TQ$ the presymplectic 2-form $\Omega=pr_1^* \omega_Q$. Let $\omega_L=-dS^*dL$ be the Poincar\'e-Cartan 2-form on $TQ$ associated to $L:TQ\rightarrow \R$ and $E_L$ its energy function. Take also the presymplectic 2-form $pr_2^*\omega_L$ on $T^*Q\oplus TQ$, and define  the functions
\begin{eqnarray*}
H&=& \langle pr_1,pr_2\rangle -  pr_2^*L \, ,\\
\bar{H}&=& \langle pr_1,pr_2\rangle -  pr_2^*E_L \, .
\end{eqnarray*}

\begin{lemma}\label{aq}
The diffeomorphism $F:T^*Q\oplus TQ\rightarrow T^*Q\oplus TQ$ is a presymplectomorphism from 
$(T^*Q\oplus TQ, \Omega)$ onto $(T^*Q\oplus TQ, \Omega+pr_2^*\omega_L)$, i.e., $F^*(\Omega+pr_2^*\omega_L)=\Omega$. Moreover, it verifies that $F^*\bar{H}=H$.
\end{lemma}
\proof{
$F$ is clearly invertible with inverse
\[
\begin{array}{cccl}
F^{-1}: & T^*Q\oplus TQ& \longrightarrow &T^*Q\oplus TQ \\
 & (\alpha, v)& \longmapsto & (\alpha+ Leg(v), v)\; .
\end{array}
\]
A direct computation shows that $H \circ F^{-1} = \bar{H}$. Moreover, in local coordinates,
\[
(F^{-1})^*(dq^A\wedge dp_A)=
dq^A\wedge \left[dp_A+d\left(\frac{\partial L}{\partial \dot q^A}\right)\right]
=dq^A\wedge dp_A+dq^A\wedge d\left(\frac{\partial L}{\partial \dot q^A}\right) \, ,
\]
which implies $F^*(\Omega+pr_2^*\omega_L)=\Omega$.
\QED
}

Denote by $j: T^*Q\times_Q M\hookrightarrow T^*Q\oplus TQ$ and $i: M\hookrightarrow TQ$ the respective canonical inclusions. Let us define $\bar{\omega}=j^*(\Omega+pr_2^*\omega_L)$. Since $pr_2 \circ j = i \circ \pi_2$, we have that
\[
\bar{\omega}=\omega+(i\circ\pi_2)^*\omega_L\; .
\]
\begin{proposition}\label{cde}
The solutions of the equations
\begin{equation}\label{zxc}
i_X\omega=dH_{W_0}\; ,
\end{equation}
 and
\begin{equation}\label{zxv}
 i_Y\bar{\omega}=d (j^*\bar{H}) \; ,
\end{equation}
are $F_{|W_0}$-related, that is, 
if $x\in T^*Q\times_Q M$ is a point where there exists a solution $Y$ of equation (\ref{zxv}) then $TF^{-1}(Y)$ is a solution of equation  (\ref{zxc}) at $F^{-1}(x)$ and, conversely, if  $X$ is a solution of equation (\ref{zxc}) at $F^{-1}(x)$ then $TF(X)$ is a solution of equation  (\ref{zxv}) at $x$.
\end{proposition}
\proof{
It readily follows from Lemma \ref{aq}.\QED}

An inmediate consequence is the following
\begin{corollary}\label{coro}
$F$ preserves the constraint submanifolds provided by the presymplectic systems $(T^*Q\times_Q M, \omega_{W_{0}}, H_{W_{0}})$ and $(T^{*}Q\times_Q M, \bar{\omega}, j^{*}\bar{H})$. That is, if 
\[
\dots \hookrightarrow W_k \dots \hookrightarrow W_1 \hookrightarrow W_0=T^*Q\times_Q M \; \; \hbox{and}
\]
\[
\dots \hookrightarrow P_k \dots \hookrightarrow P_1\hookrightarrow  P_0=T^*Q\times_Q M\,,
\]
are the sequences of submanifolds generated by the Gotay and Nester's algorithm for the first and the second presymplectic hamiltonian system, respectively, then 
\[
F_{i}=F_{|W_i}:W_i \longrightarrow P_i\,,
\]
are diffeomorphisms for all $i$.
\end{corollary}

In conclusion, Proposition \ref{cde} and Corollary \ref{coro} show that it is equivalent to solve the vakonomic hamiltonian equations (\ref{zxc}) as in Section \ref{S2} or equations (\ref{zxv}). Locally, if $(q^A(t), p_A(t), \dot{q}^a(t))$ is an integral curve of $X$ then
\[
(q^A(t), p_A-i^* \hspace{-3pt} \frac{\partial L}{\partial \dot q^A}(q^B(t), \dot{q}^b(t)),\dot{q}^a(t))
\]
is an integral curve of $Y$.

Next, we will study solutions of equations (\ref{zxv}) from a local point of view. First, notice that it is clearly equivalent to solve equations (\ref{zxv}) or the following system of equations on the ambient space
\[
\left\{
\begin{array}{l}
\displaystyle{ \left(i_{Y}( \Omega+pr_2^*\omega_L )= d \bar{H} \right)_{\big| T^*Q\times_Q M }}\,,\\
\displaystyle{ Y_{\big| T^{*}Q\times_Q M}} \in T( T^{*}Q\times_Q M)\, .
\end{array}
\right.
\]

Now, if $(q^A, \lambda_A, \dot{q}^A)$ are local coordinates on $T^*Q \oplus TQ$ and $(q^A(t), \lambda_A(t), \dot{q}^A (t))$ is an integral curve of $Y$, then, since
\[
\bar{H}= \displaystyle{ \dot{q}^A \lambda_{A} + \dot{q}^A \frac{\partial L}{\partial \dot{q}^A}}-L\, ,
\]
and
\[
\Omega+pr_2^*\omega_L=\displaystyle{ dq^A \wedge d\lambda_A+dq^A \wedge d\left(\frac{\partial L}{\partial \dot{q}^A}\right)}\; ,
\]
we have that
\begin{eqnarray*}
d \bar{H} - i_{Y }(\Omega+pr_2^*\omega_L) &=& \displaystyle{ \lambda_{A} d \dot{q}^A  + \left( \dot{q}^B \frac{\partial ^{2} L }{\partial \dot{q}^A \partial q^B}  +  \ddot{q}^B 
\frac{\partial ^{2} L }{\partial \dot{q}^A \partial \dot{q}^{B}}
-\frac{\partial L}{\partial q^A}+ \dot{\lambda}_{A}
\right) dq^{A}  = 0}\,.
\end{eqnarray*}
When restricting this formula to $ T^{*}Q\times_Q M$, we must use that $\dot{q}^{\alpha}= \Psi^{\alpha}(q^A, \dot{q}^a)$ and then, $d\dot{q}^{\alpha} = \displaystyle{\frac{ \partial \Psi^{\alpha}}{\partial q^A}}dq^A + 
\displaystyle{\frac{ \partial \Psi^{\alpha}}{\partial \dot{q}^a}}d\dot{q}^a
$. Therefore, the former system is written as
\begin{equation}\label{poi}
\left\{
\begin{array}{l}
\displaystyle{\left( \frac{d}{dt}\left( \frac{\partial L}{\partial \dot{q}^A} \right) -
\frac{\partial L}{\partial q^A} \right)_{|M} dq^A = - \dot{\lambda}_A dq^A - \lambda_{\alpha} \frac{\partial \Psi^{\alpha}}{\partial q^A} dq^{A}}\,,\\
\\
\displaystyle{\left(\lambda_{a} + \lambda_{\alpha} \frac{\partial\Psi^{\alpha}}{\partial \dot{q}^a}\right)d\dot{q}^a=0} \,.
\end{array}
\right.
\end{equation}

Moreover we have that 

\[
\dot{\lambda}_a = \frac{d}{dt}\left(-\lambda_{\alpha}\frac{\partial\Psi^{\alpha}}{\partial \dot{q}^a}\right)= -
\frac{\partial \Psi^{\alpha}}{\partial \dot{q}^a} \dot{\lambda}_{\alpha} - \lambda_{\alpha} \frac{d}{dt}
\left(\frac{\partial\Psi^{\alpha}}{\partial \dot{q}^a}\right)\; .
\]

After some computations, the system becomes

\[
\left\{
\begin{array}{l}
\displaystyle{ \frac{d}{dt}\left( \frac{\partial L}{\partial \dot{q}^A} \right) -
\frac{\partial L}{\partial q^A}   =  \dot{\lambda}_{\alpha}
\frac{\partial \Phi^{\alpha}}{\partial \dot{q}^A}
+\lambda_{\alpha}\left[\frac{d}{dt}
\left(\frac{\partial \Phi^{\alpha}}{\partial \dot{q}^A}\right)-\frac{\partial 
\Phi^{\alpha}}{\partial q^A}\right]}\; ,
  \\
\\
\Phi^{\alpha}(q, \dot{q})=0\,, \,\, 1\leq \alpha \leq m\,,
\end{array}
\right.
\]
where now $\Phi^{\alpha}(q,\dot{q})=\Psi^{\alpha}(q^A, \dot{q}^a)-\dot{q}^{\alpha}$, $1\leq \alpha\leq m$.
These last equations are the classical equations of motion for a vakonomic system or for the dynamic optimization under nonholonomic constraints (see equations (\ref{asdo})).

\subsection{Vershik-Gershkovich approach}

In \cite{VG}, Vershik and Gershkovich gave a formulation for the ``nonholonomic variational problem", i.e., the vakonomic problem, within the framework of the so-called mixed bundle, which we briefly review in the following.

If ${\cal D}:Q \longrightarrow TQ $ is a differentiable distribution along $Q$ then, the mixed bundle over $Q$ associated to ${\cal D}$ is given by ${\cal D} \oplus {\cal D}^{o}$, where ${\cal D}^{o}$ is the codistribution annihilating ${\cal D}$. This is, the fibres of ${\cal D} \oplus {\cal D}^{o} \longrightarrow Q$ are ${\cal D}_{q} \oplus {\cal D}^{o}_{q}$.

Let $\{\,\Phi^{\alpha}(q^A, \dot{q}^A)= \Psi^{\alpha}_{a}(q) \dot{q}^a-\dot{q}^{\alpha}\,, 1\leq \alpha \leq m \,\}$ be a set of independent functions whose annihilation defines the distribution ${\cal D}$ and let $\{\, \eta^{\alpha}=\Psi^{\alpha}_{a} d{q}^a- dq^{\alpha}\,, 1\leq \alpha \leq m \,\}$ be the corresponding basis of ${\cal D}^{o}$. Regarding ${\cal D} \subset TQ$ as the set of admissible velocities, Vershik and Gershkovich write the equations of motion (\ref{asdo}) for the vakonomic problem $(L,{\cal D})$ as follows
\begin{equation}\label{aqq}
\left\{
\begin{array}{l}
\displaystyle{\left( \frac{d}{dt}\left( \frac{\partial L}{\partial \dot{q}^A} \right) -
\frac{\partial L}{\partial q^A} \right)dq^A = \dot{\lambda}_{\alpha} \eta^{\alpha} + \lambda_{\alpha} (i_{\dot{q}} d \eta ^{\alpha})}\;,\\
\\
<\dot{q}, \eta^{\alpha}> = 0 \,, \,\, 1 \leq \alpha \leq m \; .
\end{array}
\right.
\end{equation}

In this particular case, we obtain that $P_1$, the first constraint submanifold for the presymplectic hamiltonian system $(T^{*}Q\times_Q M, \bar{\omega}, j^{*}\bar{H})$, is just ${\cal D}^o\oplus {\cal D}$, since
\[
\lambda_{a} + \lambda_{\alpha} \Psi_a^{\alpha}=0,\ 1\leq \alpha\leq m\; ,
\]
from equations (\ref{poi}).

If $(P_1={\cal D}^o\oplus {\cal D}, \omega_{P_1})$ is a symplectic manifold (see Proposition \ref{asdd}), then the equations of motion (\ref{aqq}) determine a unique vector field on ${\cal D} \oplus {\cal D}^o$ and the Lagrange multipliers $\lambda_{\alpha}$ are coordinates in ${\cal D}^o$ with respect to  the basis ${\eta^{\alpha}}$.

Consequently, the geometrical picture we have developed in Section \ref{S2} is equivalent to Vershik-Gershkovich approach. As said above, we have obtained a generalization of Vershik-Gershkovich formulation to the case of nonlinear constraints, just ``translating'' things from our approach by the diffeomorphism $F$. 

In the nonlinear case, under the admissibility condition, one can verify that the first constraint submanifold $P_1=F(W_1)$ can be identified with the manifold $S^*(TM^o) \times_Q M$. In fact, we have that $S^*(TM^o)$ is generated by the 1-forms 
\[
S^*d\Phi^\alpha = dq^\alpha - \frac{\partial \Psi^\alpha}{\partial \dot{q}^a}dq^a \, , \; 1 \le \alpha \le m \, . 
\]
If $(\lambda_A,q^A,\dot{q}^a) \in P_1$, then the 1-form $\lambda_A dq^A$ is a linear combination of the 1-forms $S^*d\Phi^\alpha$ in the following manner
\[
\lambda_A dq^A = \lambda_\alpha S^*d\Phi^\alpha \, .
\]

\section{Geometric approach to nonholonomic mechanics}\label{pppp}

A nonholonomic lagrangian system consists of a lagrangian $L: TQ\rightarrow \R$ subject to nonholonomic constraints defined by $m$ local functions $\Phi^{\alpha}(q^A, \dot{q}^A)$, $1\leq\alpha\leq m$.
The equations of motion for nonholonomic mechanics are derived assuming  that the constraints satisfy d'Alembert's principle, in the linear or affine case. In the nonlinear case, it does not exist an unanimous consensus about the principle to adopt \cite{Ma,Pi}. The most widely used model is the Chetaev's principle and it will be assumed in this paper. The equations of motion are then given by
\begin{equation}\label{3}
\displaystyle{\frac{d}{dt}\left(\frac{\partial L}{\partial
\dot{q}^{A}}\right) - \frac{\partial L}{\partial q^{A}} =  \lambda_{\alpha}
\frac{\partial \Phi^{\alpha}}{\partial \dot{q}^A}} \; ,
\end{equation}
together with  the algebraic equations $\Phi^{\alpha}(q^A, \dot{q}^A)=0$.
The functions $\lambda_{\alpha}$, $ 1\leq \alpha \leq m$, are some Lagrange
multipliers to be determined.
 
As in the vakonomic case, we assume the admissibility condition, so it is possible to write the constraints as
$
\dot{q}^{\alpha}=\Psi^{\alpha}(q^A, \dot{q}^a)\; ,
$
where $1\leq \alpha\leq m$, $m+1\leq a\leq n$ and $1\leq A\leq n$.

The study of nonholonomic systems in the realm of Geometric Mechanics has been an active area of research in the last years (see, for instance, \cite{MD} and references therein). The nonholonomic equations of motion can be written geometrically as
\begin{equation}
\left\{ \begin{array}{c}
(i_{\Gamma}\omega_{L} - dE_{L})_{|M} \in S^*(TM^o) \; ,\\
\Gamma_{|M} \in TM \; ,
\end{array}\right.
\label{eqc}
\end{equation}
where the subbundle $S^*(TM^o)$ of $T^*TQ$ along $M$ represents the constraint forces.

Nonholonomic mechanics also admits a nice geometrical description on the space $T^*Q\oplus TQ$ inspired in the one by Skinner and Rusk \cite{SR,SR2}. In addition, this description will be appropiate to compare the solutions of the dynamics between the vakonomic and  nonholonomic mechanics.
In the following, we will prove that equations (\ref{eqc}) are equivalent to the next ones
\begin{equation}\label{mju}
\left\{ \begin{array}{r}
\left(i_X\Omega-dH\right)_{|T^*Q\times_Q M}\in F^o\; ,\\
X_{|T^*Q\times_Q M}\in T (T^*Q\times_Q M)\; ,\\
\end{array}
\right.
\end{equation}
where
$\Omega$ is the presymplectic 2-form $\Omega=pr_1^* \omega_Q$ on $T^*Q\oplus TQ$, $H$ the hamiltonian function $H= \langle pr_1,pr_2\rangle - pr_2^*L$
and $F^o$ the subbundle of $T^*(T^*Q\oplus TQ)$ along $T^*Q\times_Q M$ defined by $F^o=pr_2^*(S^*(TM^o))$.

Indeed we have in local coordinates
\begin{eqnarray*}
\Omega&=&dq^A\wedge dp_A\; ,\\
d{H}&=&\dot{q}^Adp_A+p_Ad\dot{q}^A
-\frac{\partial L}{\partial q^A}dq^A
-\frac{\partial L}{\partial \dot{q}^A}d\dot{q}^A\; ,
\end{eqnarray*}
and $F^o$ is generated by the 1-forms
\[
\frac{\partial \Phi^{\alpha}}{\partial \dot{q}^A}dq^A
=\frac{\partial \Psi^{\alpha}}{\partial \dot{q}^a}dq^a - dq^{\alpha}
\; , \; 1\leq \alpha\leq m\; .
\]

If $\displaystyle{X=X^A\frac{\partial}{\partial q^A}+Y^A \frac{\partial}{\partial \dot{q}^A}+Z_A\frac{\partial}{\partial p_A}}$ was a solution of equations (\ref{mju}), then we would have
\begin{eqnarray}\label{1}
X^A&=&\dot{q}^A \, , \nonumber \\
Z_A&=&\frac{\partial L}{\partial q^A}+\lambda_{\alpha}\frac{\partial \Phi^{\alpha}}{\partial \dot{q}^A} \, ,
\end{eqnarray}
\noindent along with the constraints
\begin{eqnarray}\label{2}
p_A-\frac{\partial L}{\partial \dot{q}^A}&=&0\; , \nonumber \\
\Phi^{\alpha}(q^A, \dot{q}^A)&=&0\; .
\end{eqnarray}

Observe that these constraints determine the submanifold $\tilde{M}$ of $T^*Q\times_Q M$. The submanifold $\tilde{M}$ is diffeomorphic to $M$ since
\[
\begin{array}{rcl}
M&\longrightarrow&\tilde{M}\\
m&\longmapsto&(Leg_L(m), m) \, ,
\end{array}
\]
is a diffeomorphism. $\tilde{M}$ is the first constraint submanifold provided by the constraint algorithm applied to equations (\ref{mju}). This algorithm will lead to a final constraint submanifold on which there exists a well-defined dynamics, at least in case the given problem is consistent (see \cite{MD}). Obviously, equations (\ref{1}) and (\ref{2}) are equivalent to the nonholonomic equations of motion (\ref{3}).

In terms of the $\Psi^{\alpha}$'s the above equations are written as 
\begin{eqnarray*}
X^A&=&\dot{q}^A\; ,\\
Z_a&=&\frac{\partial L}{\partial q^a}+\lambda_{\alpha}\frac{\partial \Psi^{\alpha}}{\partial \dot{q}^a}\; ,\\
Z_{\beta}&=&\frac{\partial L}{\partial q^{\beta}}-\lambda_{\beta}\; ,
\end{eqnarray*}
along with the constraints
\begin{eqnarray}\label{oiu}
p_A-\frac{\partial L}{\partial \dot{q}^A}&=&0\; , \nonumber \\
\dot{q}^{\alpha}-\Psi^{\alpha}(q^A, \dot{q}^a)&=&0\; .
\end{eqnarray}

Therefore, a solution $X$ of (\ref{mju}) is of the form
\begin{eqnarray*}
X &=& \dot{q}^a \left( \frac{\partial }{\partial q^a} + \frac{\partial \Psi^\alpha}{\partial q^a} \frac{\partial}{\partial \dot{q}^\alpha} + \left( \frac{\partial^2 L}{\partial \dot{q}^A \partial q^a} + \frac{\partial \Psi^\alpha}{\partial q^a} \frac{\partial^2 L}{\partial \dot{q}^A \partial \dot{q}^\alpha} \right) \frac{\partial}{\partial p_A} \right) \\
&&+ \Psi^\gamma \left( \frac{\partial }{\partial q^\gamma} + \frac{\partial \Psi^\alpha}{\partial q^\gamma} \frac{\partial}{\partial \dot{q}^\alpha} + \left( \frac{\partial^2 L}{\partial \dot{q}^A \partial q^\gamma} + \frac{\partial \Psi^\alpha}{\partial q^\gamma} \frac{\partial^2 L}{\partial \dot{q}^A \partial \dot{q}^\alpha} \right) \frac{\partial}{\partial p_A} \right) \\
&& +Y^a \left( \frac{\partial }{\partial \dot{q}^a} + \frac{\partial \Psi^\alpha}{\partial \dot{q}^a} \frac{\partial}{\partial \dot{q}^\alpha} + \left( \frac{\partial^2 L}{\partial \dot{q}^A \partial \dot{q}^a} + \frac{\partial \Psi^\alpha}{\partial \dot{q}^a} \frac{\partial^2 L}{\partial \dot{q}^A \partial \dot{q}^\alpha} \right) \frac{\partial}{\partial p_A} \right) \, .
\end{eqnarray*}

Under the regularity assumption, which means that the matrix
\begin{equation}\label{6}
\tilde{\cal C}_{ab}=\frac{\partial^2 \tilde{L}}{\partial \dot{q}^a\partial \dot{q}^b}-
i^* \hspace{-3pt}\left(\frac{\partial {L}}{\partial \dot{q}^{\alpha}}\right)
\frac{\partial^2 \Psi^{\alpha}}{\partial \dot{q}^a \partial \dot{q}^b} \, ,
\end{equation}
is invertible (see \cite{WSF}), there is an unique solution of the dynamics on $\tilde{M}$. In particular, after some computations, we obtain
\begin{eqnarray}\label{ytr1}
Y^a&=&-\tilde{\cal C}^{ab}\left[
\dot{q}^A\frac{\partial^2 \tilde{L}}{\partial q^A\partial\dot{q}^b}-
\dot{q}^A 
i^* \hspace{-3pt} \left(\frac{\partial L}{\partial\dot{q}^{\alpha}}\right)
\frac{\partial^2 \Psi^{\alpha}}{\partial q^A \partial \dot{q}^b}\right. \nonumber \\
&&\left.-\frac{\partial\tilde{L}}{\partial q^b}+i^* \hspace{-3pt} \left(\frac{\partial L}
{\partial q^{\alpha}}\right)\left(
\frac{\partial \Psi^{\alpha}}{\partial {q}^b}-
\frac{\partial \Psi^{\alpha}}{\partial \dot{q}^b}\right)\right]\; ,
\end{eqnarray}
where $i:M\rightarrow TQ$ is the canonical inclusion and $\tilde{\cal C}^{ab}$ the inverse matrix of $\tilde{\cal C}_{ab}$.

Taking coordinates $(q^A, \dot{q}^a)$ on $\tilde{M}$, the equations of motion for a nonholonomic system will be 
\begin{equation}\label{4}
\left\{
\begin{array}{rcl}
\dot{q}^{\alpha}&=&\Psi^{\alpha}(q^A, \dot{q}^a) \, , \\
\ddot{q}^a&=&\displaystyle{-{\cal C}^{ab}\left[
\dot{q}^A\frac{\partial^2 \tilde{L}}{\partial q^A\partial\dot{q}^b}-
\dot{q}^A i^* \hspace{-3pt} \left(\frac{\partial L}{\partial\dot{q}^{\alpha}}\right)
\frac{\partial^2 \Psi^{\alpha}}{\partial q^A \partial \dot{q}^b} \right.}  \\
&& -\displaystyle{ \left. \frac{\partial\tilde{L}}{\partial q^b}+i^* \hspace{-3pt} \left(\frac{\partial L}
{\partial q^{\alpha}}\right)\left(
\frac{\partial \Psi^{\alpha}}{\partial {q}^b}-
\frac{\partial \Psi^{\alpha}}{\partial \dot{q}^b}\right)\right]}
\end{array}
\right.
\end{equation}
\noindent Compare them with equations (\ref{ytr}).

\section{Vakonomic and nonholonomic mechanics: Equivalence of dynamics}

In this section, we shall investigate the relation between vakonomic and nonholonomic dynamics.

Consider a physical system with lagrangian $L: TQ \rightarrow
\R$ and constraint submanifold $M \subset TQ$. Let us assume that the
vakonomic problem lives in the first constraint submanifold, $W_1$, and
that the nonholonomic one lives in $\tilde{M}$ (this will be the case
if the constraints are linear and the admissibility and compatibility
conditions are satisfied). As a consequence, we have well defined vector
fields $X_{vk}$ on $W_1$ and $X_{nh}$ on $\tilde{M}$.

It is clear that the mapping $(\pi_2)_{|W_1}: W_1\rightarrow M$ is a surjective submersion and that we can define  the mapping $\Upsilon: W_1\rightarrow \tilde{M}$ as
\[
\begin{array}{rccl}
\Upsilon:& W_1 & \longrightarrow & \tilde{M} \\
& (\alpha, v)&\longmapsto&(Leg_L(v), v)
\end{array}
\]
In coordinates, $\Upsilon$ reads as $\Upsilon(q^A,\dot{q}^a,p_{\alpha})= (q^A,\dot{q}^a)$.

Our aim is to know whether, given a nonholonomic solution, we can find
initial conditions in the vakonomic Lagrange multipliers, $p_{\alpha}$,
so that the curve can also be seen as a vakonomic solution. In order to
capture the common solutions to both problems, we have developed the
following algorithm. It is inspired in the idea of the $\Upsilon$-relation of $X_{vk}$ and $X_{nh}$ and the constraint algorithm developed by O. Krupkov\'a \cite{K}. If both fields were $\Upsilon$-related, then the projection to $\tilde{M}$ of all the vakonomic solutions would be nonholonomic. So, selecting the points in which they are related, we are picking up all the possible good candidates. We write $W_1=S_0$ and define
\[
S_1=\{ w \in S_0 \, | \; T_w \Upsilon (X_{vk}(w)) = X_{nh}(\Upsilon (w)) \} \, .
\]

In general $S_1$ is not a submanifold. If $S_1=\emptyset$, there is no relation between the vakonomic and nonholonomic dynamics.

If $S_1\not=\emptyset$, we apply the following algorithm:

\begin{itemize}
\item {\bf Step 1}: For any $w\in S_1$, consider $C_{(w)}=\cup_{i}C_{(w)i}$, the union of all connected submanifolds $C_{(w)i}$ of maximal dimension lying in $S_1$, contained in a neighbourhood $U$ of $w$ and passing through $w$ (maximal dimension means that if $N$ is a connected submanifold lying in $S_1\cap U$ passing through $w$ and $C_{(w)i}\subseteq N$ then $C_{(w)i}=N$).

Suppose that $C_{(w)}\not= \{w\}$. For each $i$ we consider the subset of $C_{(w)i}$
\[
\tilde{C}_{(w)i}=\{ v \in C_{(w)i} \, | \; X_{vk}(v)\in T_vC_{(w)i}\}  \, .
\]
If $\tilde{C}_{(w)i}={C}_{(w)i}$ then we call the submanifold $C_{(w)i}$ a {\bf final constraint submanifold at $w$}. If $\tilde{C}_{(w)i}=\emptyset$, we exclude $C_{(w)i}$ from the bunch $C_{(w)}$. If $\emptyset \subsetneq \tilde{C}_{(w)i}\subsetneq C_{(w)i}$, then we proceed to the next step.

\item {\bf Step 2}: Repeat the Step 1 with  $\tilde{C}_{(w)i}$ instead of $S_1$.

\end{itemize}

After sufficient steps of this algorithm we obtain a bunch of final constraint submanifolds at $w$ or we find that there is no final constraint submanifold passing through $w$. Collecting all the points where there exists a bunch of final constraint submanifolds we obtain the subset where there is equivalence between vakonomic and noholonomic dynamics. 
 
Suppose that the constraints $\Phi^{\alpha}$, $1\leq \alpha\leq m$, are linear on the velocities so we can write them as
\[
\dot{q}^{\alpha}=\Psi^{\alpha}_a(q)\dot{q}^a\; . 
\]
In such a case, the matrices ${\cal C}$ and $\tilde{\cal C}$ defined in (\ref{5}) and (\ref{6}), respectively, are the same (even for constraints affine on the velocities).

\begin{proposition}
$S_1$ is locally chararacterized by the vanishing of the $n-m$ constraints functions on $W_1$
\begin{equation}\label{8}
g_b=\dot{q}^a\left(p_{\alpha}-i^* \hspace{-3pt} \frac{\partial L}{\partial \dot{q}^{\alpha}}\right)\left[
\frac{\partial \Psi^{\alpha}_b}{\partial q^a}
-\frac{\partial \Psi^{\alpha}_a}{\partial q^b}
+\Psi^{\beta}_{a}\frac{\partial \Psi^{\alpha}_b}{\partial q^{\beta}}
-\Psi^{\beta}_{b}\frac{\partial \Psi^{\alpha}_a}{\partial q^{\beta}}
\right], \quad m+1\leq b\leq n\; .
\end{equation}
\end{proposition}
\proof{ The comparison between the vector fields $X_{vk}$ and $X_{nh}$ consists of taking the difference between $\ddot{q}^a$'s in the expressions (\ref{ytr}) and (\ref{4}) and equating the result to zero. \QED}

Consider the local projection $\rho(q^a,q^{\alpha})=(q^{\alpha})$ and the connection $\Gamma$ on $\rho$ such that the horizontal distribution ${\cal H}$ is given by prescribing its annihilator to be
\[
{\cal H}^o=\langle dq^{\alpha}-\Psi^{\alpha}_a dq^a, 1\leq \alpha\leq m\rangle\; .
\]
Then the curvature $R$ of this connection (see \cite{LR}) is given by
\[
R(\frac{\partial}{\partial q^{a}}, \frac{\partial}{\partial
q^{b}}) = R^{\alpha}_{ab} \frac{\partial}{\partial q^{\alpha}} \; ,
\]
where
\[
R^{\alpha}_{ab} = 
\frac{\partial \Psi^{\alpha}_b}{\partial q^a}
-\frac{\partial \Psi^{\alpha}_a}{\partial q^b}
+\Psi^{\beta}_{a}\frac{\partial \Psi^{\alpha}_b}{\partial q^{\beta}}
-\Psi^{\beta}_{b}\frac{\partial \Psi^{\alpha}_a}{\partial q^{\beta}}\; .
\]
We say that $\Gamma$ is flat if the curvature $R$ vanishes identically. The tensor $R$ measures the lack of integrability of the horizontal distribution ${\cal H}$, which in our case is the constraint manifold.

Then, we can write the constraints determining $S_1$ as
\[
g_b=\dot{q}^a\left(p_{\alpha}-i^* \hspace{-3pt} \frac{\partial L}{\partial \dot{q}^{\alpha}}\right) R^{\alpha}_{ab}, \quad m+1\leq b\leq n\; .
\]

>From this expression we obtain that if the constraints are holonomic, then $R=0$ and the final constraint submanifold is equal to $S_0=W_1$. Therefore, every nonholonomic solution is also a vakonomic solution. Indeed, equations (\ref{asx}) will read as
\begin{equation}\label{asx1}
\left\{
\begin{array}{l}
\dot{q}^{\alpha}=\Psi^{\alpha}_a\dot{q}^a\; ,\\
\displaystyle{\dot{p}_{\alpha}= \frac{\partial \tilde{L}}{\partial q^{\alpha}}-p_{\beta}\frac{\partial \Psi^{\beta}_a}{\partial q^{\alpha}}\dot{q}^a}\; ,\\
\displaystyle{\frac{d}{dt}
\left(\frac{\partial \tilde{L}}{\partial \dot{q}^a}\right)
-\frac{\partial \tilde{L}}{\partial q^{a}}=\Psi^{\alpha}_a\frac{\partial \tilde{L}}{\partial q^{a}}} \; .
\end{array}
\right.
\end{equation}
The first and the third set of equations determine the trajectory in $M$. The Lagrange multipliers $p_{\alpha}$ are determined by the second set of equations once we know the solution in $M$. This is the typical behavior of the holonomic case \cite{LeMu}.

But, in general, for linear constraints, the first constraint subset in the algorithm is determined by
\[
S_1 = \{ g_b=0,\ m+1\leq b\leq n  \} \, , 
\]
where $g_b(q^A,\dot{q}^a,p_\alpha) = \displaystyle{\dot{q}^a R^\alpha_{ab}(q)(p_\alpha - \frac{\partial L}{\partial \dot{q}^\alpha})}$. Note that $S_1$ will not be a submanifold, because $0$ is not a regular value of the functions $g_b$, $b=m+1,...,n$. 
Anyway, the geometric context we have developed can be very useful to tackle the problem of the comparison of the two methods.

\begin{proposition}\label{9}
If $c(t)=(q^A(t))$ is a solution of the free problem 
which verifies all the constraints, i.e,
\[
\dot{q}^{\alpha}(t)=\Psi^{\alpha}_a(q(t))\dot{q}^a(t)\, , \; 1 \leq \alpha \leq m \, ,
\]
then $c(t)$ is a solution of the nonholonomic and vakonomic problems simultaneously.
\end{proposition}

\proof{ Let us consider the submanifold $S=\{ p_\alpha = i^* \hspace{-3pt} \displaystyle{
\frac{\partial L}{\partial \dot{q}^\alpha}} \}$, which is contained in
$S_1$. A natural question is whether the vakonomic vector field
will be tangent to $S$, that is, $X_{vk} \in TS$.
>From equations (\ref{asx}), we have  along any integral curve of the vakonomic vector field
\[ 
X_{vk} \in TS 
\Longleftrightarrow
\frac{d}{dt}\left(p_\alpha - i^* \hspace{-3pt} \frac{\partial L}{\partial \dot{q}^\alpha}\right) = 0
\Longleftrightarrow 
\dot{p}_\alpha = \dot{q}^A \frac{\partial^2 L}{\partial q^A \partial \dot{q}^\alpha} + \ddot{q}^a \frac{\partial^2
L}{\partial \dot{q}^a \partial \dot{q}^\alpha} \, .
\] 
On $S$, we have that
\[
\dot{p}_{\alpha} = \displaystyle{\frac{\partial \tilde{L}}{\partial q^\alpha} - p_\beta \frac{\partial \Psi^\beta}{\partial q^\alpha}= \frac{\partial \tilde{L}}{\partial q^\alpha} - \frac{\partial L}{\partial \dot{q}^\beta} \frac{\partial \Psi^\beta}{\partial q^\alpha} = \frac{\partial L}{\partial q^\alpha}} \, .
\]
Then the above condition can be rewritten as 
\[ 
\frac{\partial L}{\partial q^\alpha} =
\frac{d}{dt} \left( \frac{\partial L}{\partial \dot{q}^\alpha} \right)
\; ,
\] 
that with the third set of equations in (\ref{asx}) are precisely the Euler-Lagrange equations. Then, we have proved that  $c(t)$ is a free solution and satisfies the constraints if and only if
\[
(q^A(t), i^* \hspace{-3pt} \frac{\partial L}{\partial \dot{q}^{\alpha}}, \dot{q}^a(t) )\; ,
\]
is a solution of the vakonomic equations (\ref{asx}). Since the constraints $g_b=0$ are automatically satisfied for all the points in $S$ we deduce that $c(t)$ is also a solution of the nonholonomic problem. \QED
}

\begin{remark}
{\rm 
As a consequence of Proposition \ref{9} we obtain that
if $g$ is a riemannian metric on $Q$ with kinetic energy $L=\frac{1}{2}g$ and we assume that we are given 
a distribution ${\cal D}$ on $Q$ which is geodesically invariant with respect to the Levi-Civita connection $\nabla^g$, then all the nonholonomic solutions can be seen as vakonomic ones. In fact, they all are solutions of the free problem. This last result was first stated in \cite{Fa} (Theorem 3.2) with additional hypothesis on the nature of the metric $g$ and the integrability of ${\cal D}^{\perp_g}$ which are not essential, as we have seen.
}
\end{remark}

\begin{remark}
{\rm 
Let $\Theta: G \times Q \longrightarrow Q$ be a free and proper action on $Q$. Then $\pi: Q \longrightarrow Q/G$ is a principal $G$-bundle. Assume that the lagrangian $L:TQ \longrightarrow \R$ is $G$-invariant and is subject to equivariant affine constraints, $M$, such that its linear part ${\cal D}$ is the horizontal distribution of a principal connection $\gamma$ on $\pi: Q \longrightarrow Q/G$. Then, we have the following result, which is an adaptation of Theorem 3.1 in \cite{Fa} to our geometric description of vakonomic and nonholonomic mechanics.

\begin{proposition}
Assume that the admissibility and compatibility conditions hold. Then, the following are equivalent:
\begin{enumerate}
\item the nonholonomic solution $(q^A(t),\dot{q}^a(t)) \in \tilde{M}$ verifies the condition $g_b(q^A(t),\dot{q}^a(t),p_0)=0$ for some $p_0$, $m+1 \le b \le n$.
\item the curve $(q^A(t),\dot{q}^a(t),p_0) \in W_1$ is a vakonomic solution.
\end{enumerate}
\end{proposition}
}
\end{remark}

\begin{example}[Rolling penny \cite{BlCr}]
{\rm Consider a vertical penny constrained to roll without slipping on an horizontal plane and free to rotate about its vertical axis. Let $x$ and $y$ denote the position of contact of the disk in the plane. The remaining variables are $\theta$ denoting the orientation of a chosen material point $P$ with respect to the vertical and $\phi$ the heading angle of the penny. The configuration space is then $Q=\R^2 \times S^1 \times S^1$. The lagrangian defining the dynamical problem may be written as 
\[
L=\frac{1}{2}\left(\dot{x}^2+\dot{y}^2+\dot{\theta}^2+\dot{\phi}^2
\right)
\]
and the constraints are given by
\begin{eqnarray*}
\dot{x}&=&\dot{\theta}\cos\phi\; ,\\
\dot{y}&=&\dot{\theta}\sin\phi\; .
\end{eqnarray*}
For simplicity, we assume that the mass $m$, the moments of inertia $I$, $J$ and the radius of the penny $R$ are $1$.

We have well defined vector fields
\begin{eqnarray*}
X_{vk} &=& \dot{x} \frac{\partial}{\partial x} + \dot{y} \frac{\partial}{\partial y} + \dot{\theta} \frac{\partial}{\partial \theta} + \dot{\phi} \frac{\partial}{\partial \phi} + \frac{1}{2} \dot{\phi}\left(p_x \sin \phi \cos \phi -p_y\cos^2 \phi - 2 \dot{\theta}\sin \phi\right) \frac{\partial}{\partial \dot{x}}\\
&& + \frac{1}{2} \dot{\phi}\left( p_x \sin^2 \phi - p_y \cos \phi \sin \phi + 2 \dot{\theta}\cos \phi \right) \frac{\partial}{\partial \dot{y}} + \frac{1}{2}\dot{\phi} \left( p_x\sin \phi -p_y\cos \phi \right) \frac{\partial}{\partial \dot{\theta}} \\
&& - \dot{\theta} \left( p_x \sin \phi - p_y \cos \phi \right) \frac{\partial}{\partial \dot{\phi}} + \dot{\phi}\left( p_x \sin \phi \cos \phi - p_y\cos^2 \phi - 2 \dot{\theta}\sin \phi\right) \frac{\partial}{\partial p_x} \\
&& + \dot{\phi}\left(p_x \sin^2 \phi - p_y \cos \phi \sin \phi + 2 \dot{\theta}\cos \phi\right) \frac{\partial}{\partial p_y}
\end{eqnarray*}
on $W_1$ and
\begin{eqnarray*}
X_{nh} = \dot{x} \frac{\partial}{\partial x} + \dot{y} \frac{\partial}{\partial y} + \dot{\theta} \frac{\partial}{\partial \theta} + \dot{\phi} \frac{\partial}{\partial \phi} - \dot{\phi} \, \dot{\theta}\sin \phi \frac{\partial}{\partial \dot{x}} + \dot{\phi} \, \dot{\theta}\cos \phi \frac{\partial}{\partial \dot{y}}
\end{eqnarray*}
on $\tilde{M}$, respectively. Thus, we find that $S_1$ is determined by the vanishing of
\begin{eqnarray*}
g_1&=&\dot{\phi}\left( p_x \sin\phi - p_y \cos\phi \right) \; ,\\
g_2&=&-\dot{\theta}\left( p_x \sin\phi - p_y \cos\phi \right) \; .
\end{eqnarray*}
$S_1$ has two connected components
\begin{eqnarray*}
C_1 &=& \{ w \in W_1 \, | \; p_x \sin\phi - p_y \cos\phi = 0 \} \, , \\
C_2 &=& \{ w \in W_1 \, | \; \dot{\theta} = 0 \, , \; \dot{\phi} = 0 \} \, .
\end{eqnarray*}
Applying the algorithm, we obtain that $\tilde{C}_1 = C_{11} \cup C_{12}$ and $\tilde{C}_2 = C_2$, where
\begin{eqnarray*}
C_{11} &=& \{ w \in C_1 \, | \; \dot{\phi} = 0 \} \, , \\
C_{12} &=& \{ w \in C_1 \, | \; 2\dot{\theta} = p_x \cos \phi + p_y \sin \phi \} \, .
\end{eqnarray*}
One step more yields that $\tilde{C}_{11}=C_{11}$ and $\tilde{C}_{12}=C_{12}$, so they are also final constraint submanifolds. The nonholonomic solutions that fall into $C_{11}$ are motions of the penny along a straight line in the horizontal plane. The nonholonomic solutions in $C_2$  are stationary positions. However, any nonholonomic solution can be seen as a vakonomic one contained in the final constraint manifold $C_{12}$, with Lagrange multipliers $p_x= 2 \dot{\theta} \cos \phi$ and $p_y = 2 \dot{\theta} \sin \phi$. In terms of the extended lagrangian formalism mentioned in Remark \ref{exten}, we have the following Lagrange multipliers
\begin{eqnarray*}
\lambda_x &=& \frac{\partial L}{\partial x} - p_x = \dot{x} - p_x = - \dot{\theta} \cos \phi \, , \\
\lambda_y &=& \frac{\partial L}{\partial y} - p_y = \dot{y} - p_y = - \dot{\theta} \sin \phi \, , 
\end{eqnarray*}
which is just the result of Bloch and Crouch \cite{BlCr}.
}
\end{example}

\begin{example}[Constrained particle \cite{Ro}]
{\rm We will discuss here an instructive example which has been extensively treated in the literature of nonholonomic mechanics. Consider a particle of unit mass moving in space, $Q=\R^3$, subject to the constraint
\[
\Phi=\dot{z}-y\dot{x} = 0\, .
\]
The lagrangian is
\[
L=\frac{1}{2} \left( \dot{x}^2 + \dot{y}^2 + \dot{z}^2 \right) \, .
\]
The constraint is linear and the compatibility and admissibility conditions are satisfied, so we have well defined vector fields
\[
X_{vk}= \dot{x}\frac{\partial }{\partial x} + \dot{y} \left( \frac{\partial }{\partial y} + \dot{x} \frac{\partial }{\partial \dot{z}} \right) + \dot{z} \frac{\partial }{\partial z} + \frac{\dot{y}(p_z-\dot{z})-y\dot{y}\dot{x}}{1+y^2} \left( \frac{\partial }{\partial \dot{x}} +y \frac{\partial }{\partial \dot{z}} \right) - \dot{x}(p_z - \dot{z})\frac{\partial}{\partial \dot{y}}
\]
on $W_1$ and
\[
X_{nh} = \dot{x}\frac{\partial }{\partial x} + \dot{y} \left( \frac{\partial }{\partial y} + \dot{x} \frac{\partial }{\partial \dot{z}} \right) + \dot{z} \frac{\partial }{\partial z} - \frac{y\dot{y}\dot{x}}{1+y^2} \left( \frac{\partial }{\partial \dot{x}} +y \frac{\partial }{\partial \dot{z}} \right)
\]
on $\tilde{M}$, respectively.
Comparing them, we find that $S_1$ is determined by the vanishing of
\begin{eqnarray*}
g_1 &=& \dot{x}(p_z - \dot{z}) \, ,\\
g_2 &=& -\dot{y}(p_z - \dot{z}) \, .
\end{eqnarray*}
Consequently, $S_1$ has two connected components
\begin{eqnarray*}
C_1 &=& \{ w \in W_1 \, | \; p_z - \dot{z} = 0 \} \, ,\\
C_2 &=& \{ w \in W_1 \, | \; \dot{x} = 0 \, , \; \dot{y} = 0 \} \, .
\end{eqnarray*}
Applying the algorithm, we obtain that $\tilde{C}_1=C_{11} \cup C_{12}$, where
\begin{eqnarray*}
C_{11} &=& \{ w \in C_1 \, | \; \dot{y} = 0 \} \, ,\\
C_{12} &=& \{ w \in C_1 \, | \; \dot{x} = 0 \} \, .
\end{eqnarray*}
On the other hand, $\tilde{C}_2=C_2$, so $C_2$ is a final constraint submanifold. Another step of the algorithm yields
\[
\tilde{C}_{11} = C_{11} \, , \quad \tilde{C}_{12} = C_{12} \, ,
\]
so they both are also final constraint submanifolds.

Therefore, the nonholonomic solutions that can be seen as vakonomic ones are the ones which belong to
\[
\begin{array}{rl}
C_{11}: & (\dot{x}_0 t + x_0,y_0,\dot{z}_0 t + z_0) \, , \; \hbox{where} \; \dot{z}_0 = y_0 \dot{x}_0 \, ,\\
C_{12}: & (x_0,\dot{y}_0 t + y_0, z_0) \, , \\
C_2: & (x_0,y_0,z_0) \, ,
\end{array}
\]
that is, stationary or free motions in $M$. Observe that there are plenty of nonholonomic solutions that can not be seen as vakonomic ones.
}
\end{example}

\begin{example}[Ball on a rotating table \cite{LeMu}]
{\rm Applying the algorithm to this example, one can obtain the same result found in \cite{LeMu}. The configuration space is $Q=\R^2 \times SO(3)$ with coordinates $(x,y,R)$. We denote the spatial angular velocity by $\xi \in \R^3$, where $\hat{\xi} = \dot{R}R^T$. The lagrangian is
\[
L= \frac{1}{2} I \left( (\xi^1)^2 + (\xi^2)^2 + (\xi^3)^2 \right) + \frac{1}{2} m \left( \dot{x}^2 + \dot{y}^2 \right) \, ,
\]
where $I$ and $m$ are the inertia and mass of the ball, respectively. The constraints are
\begin{eqnarray*}
\dot{x} &=& r \xi^2 - \Omega y \, , \\
\dot{y} &=& -r \xi^1 + \Omega x \, ,
\end{eqnarray*}
where $r$ is the radius of the ball and $\Omega$ is the angular velocity of the table.

Applying the algorithm, one finds the following final constraint submanifolds
\begin{eqnarray*}
{C_f}_1 &=& \{ w \in W_1 \, | \; \dot{x} = \dot{y} = p_x = p_y = 0 \} \, , \\
{C_f}_2 &=& \{ w \in W_1 \, | \; \xi^3 = \Omega \} \, .
\end{eqnarray*}
As proved in \cite{LeMu}, there are nonholonomic solutions that can not be seen as vakonomic ones.
}
\end{example}

\section*{Acknowledgements}

This work was partially supported by grant DGICYT (Spain) PB97-1257. 
J. Cort\'es and S. Mart{\'\i}nez wish to thank the Spanish Ministerio de 
Educaci\'on y Cultura for FPU and FPI grants, respectively. They and D. Mart{\'\i}n de Diego would like to thank the Department of Mathematical Physics and Astronomy of the University of Ghent for its kind hospitality. We would like to thank F. Cantrijn and E. Mart{\'\i}nez for several helpful conversations.

{\parindent 0.6cm

{\sc Jorge Cort{\'e}s \dag, Manuel de Le\'on \ddag\  and Sonia Mart{\'\i}nez \S} 

{\it Laboratory of Dynamical Systems, Mechanics and Control, 
Instituto de Matem\'aticas y F{\'\i}sica Fundamental, 
Consejo Superior de Investigaciones Cient{\'\i}ficas, 
Serrano 123, 28006 Madrid, SPAIN\\
\dag e-mail: j.cortes@imaff.cfmac.csic.es\quad
\ddag e-mail: mdeleon@imaff.cfmac.csic.es\\
\S e-mail: s.martinez@imaff.cfmac.csic.es\\

{\sc David Mart{\'\i}n de Diego}

{\it Departamento de Econom{\'\i}a Aplicada (Matem\'aticas), 
Facultad de CC. Econ\'omicas y Empresariales, 
Universidad de Valladolid,
Avda. Valle Esgueva 6, 47011 Valladolid, SPAIN\\
e-mail: dmartin@esgueva.eco.uva.es}
    
\end{document}